\documentclass[11pt]{article}
\usepackage{amsmath,amssymb, color}
\usepackage{graphicx} 
\usepackage{underscore}
\usepackage{amsthm,epsfig,mathrsfs}
\usepackage{amsbsy,multirow,caption}
\usepackage{hyperref}
\usepackage{bookmark}

\newtheorem{theorem}{Theorem}

\newtheorem{ruul}[theorem]{Rule}

\newtheorem{const}{Construction}
\newtheorem{ex}{Example}

\newcounter{bulletlist}

\newcounter{enumlist0}
\newenvironment{enumlist0}{
  \begin{list}
  {(\arabic{enumlist0})}
  {\usecounter{enumlist0}
  \setlength{\leftmargin}{\parindent}
  \setlength{\labelsep}{0.6em}
  \setlength{\labelwidth}{1em}
  \setlength{\topsep}{1ex}
  \setlength{\rightmargin}{0em}
  \setlength{\itemsep}{0.8ex}
  \setlength{\parsep}{0em}
  \setlength{\itemindent}{0em} }}
  {\end{list}}

\newcounter{enumlist1}
\newenvironment{enumlist1}{
  \begin{list}
  {\arabic{enumlist1})}
  {\usecounter{enumlist1}
  \setlength{\leftmargin}{0em}
  \setlength{\labelsep}{0.6em}
  \setlength{\labelwidth}{-0.5em}
  \setlength{\topsep}{1ex}
  \setlength{\rightmargin}{0em}
  \setlength{\itemsep}{0.8ex}
  \setlength{\parsep}{0em}
  \setlength{\itemindent}{0em} }}
  {\end{list}}

\newlength\enumsep
\newcounter{enumlist}
\newenvironment{enumlist}{
  \begin{list}
  {(\alph{enumlist})}
  {\usecounter{enumlist}
  \setlength{\leftmargin}{\enumsep}
  \setlength{\labelsep}{0.6em}
  \setlength{\labelwidth}{1em}
  \setlength{\topsep}{1ex}
  \setlength{\rightmargin}{0em}
  \setlength{\itemsep}{0.5ex}
  \setlength{\parsep}{0em}
  \setlength{\itemindent}{0em} }}
  {\end{list}}

\newcommand{\R}{\mathbb{R}}

\newcommand{\circc}{\hspace{-2pt}\circ\hspace{-2pt}}

\begin{document}

\begin{center}
{\bfseries \Large  Lorentz-Conformal Transformations in the Plane} \\
\end{center}

\vspace{0.2ex}

\begin{center}
{\large Barbara A. Shipman, Patrick D. Shipman, and Stephen P. Shipman} \\
\vspace{1ex}
{\itshape Departments of Mathematics at the University of Texas at Arlington, Colorado State University in Fort Collins, and Louisiana State University in Baton Rouge}
\end{center}

\vspace{2ex}

\begin{abstract}
\noindent While conformal transformations of the plane preserve Laplace's equation, Lorentz-conformal mappings preserve the wave equation.   We discover how simple geometric objects, such as quadrilaterals and pairs of crossing curves, are transformed under nonlinear Lorentz-conformal mappings. 
Squares are transformed into curvilinear quadrilaterals where three sides determine the fourth by a geometric ``rectangle rule," which can be expressed also by functional formulas.  There is an explicit functional degree of freedom in choosing the mapping taking the square to a given quadrilateral.
We characterize classes of Lorentz-conformal maps by their symmetries under subgroups of the dihedral group of order eight.  Unfoldings of non-invertible mappings into invertible ones are reflected in a change of the symmetry group. 
The questions are simple; but the answers are not obvious, yet have beautiful geometric, algebraic, and functional descriptions and proofs.  This is due to the very simple form of nonlinear Lorentz-conformal transformations in dimension $1\!+\!1$, provided by characteristic coordinates.  

 \end{abstract}

\vspace{2ex}
\noindent
\begin{mbox}
{\bf Key words:}  Lorentz transformation; Lorentz-conformal; Lorentz invariance; wave equation
\end{mbox}
\vspace{2ex}

\hrule
\vspace{1.1ex}

\section{Lorentz-conformal maps}  

Which transformations $(x,y) \mapsto (u,v)$ preserve the wave equation,
so that $f_{xx}=f_{yy} \implies f_{uu}=f_{vv}$ ?
Characteristic coordinates
\begin{equation}\label{coordinates}
  \renewcommand{\arraystretch}{1.1}
\left.
  \begin{array}{ll}
    X=x+y \\
    Y=-x+y
  \end{array}
\right.
\qquad
  \renewcommand{\arraystretch}{1.1}
\left.
  \begin{array}{ll}
    U=u+v \\
    V=-u+v
  \end{array}
\right.
\end{equation}
allow for a simple answer.
They reduce the wave equation $f_{xx}=f_{yy}$ to  $f_{XY}=0$, and the
solution has the general form $f(X,Y)=f_1(X) + f_2(Y)$.  The transformations $(U,V)=\alpha(X,Y)$ that take these solutions into solutions $f(U,V)=g_1(U)+g_2(V)$ of $f_{UV}=0$ are
those that decouple the characteristic coordinates.  This means that either
\begin{eqnarray*}
  && (U,V) = \alpha(X,Y) \doteq (h(X),k(Y)) \quad \text{or} \\
  && (U,V) = \alpha(X,Y) \doteq (k(Y),h(X))
\end{eqnarray*}
for some functions $h$ and $k$.

\begin{figure}[h]
\centerline{\includegraphics[width=5.9in]{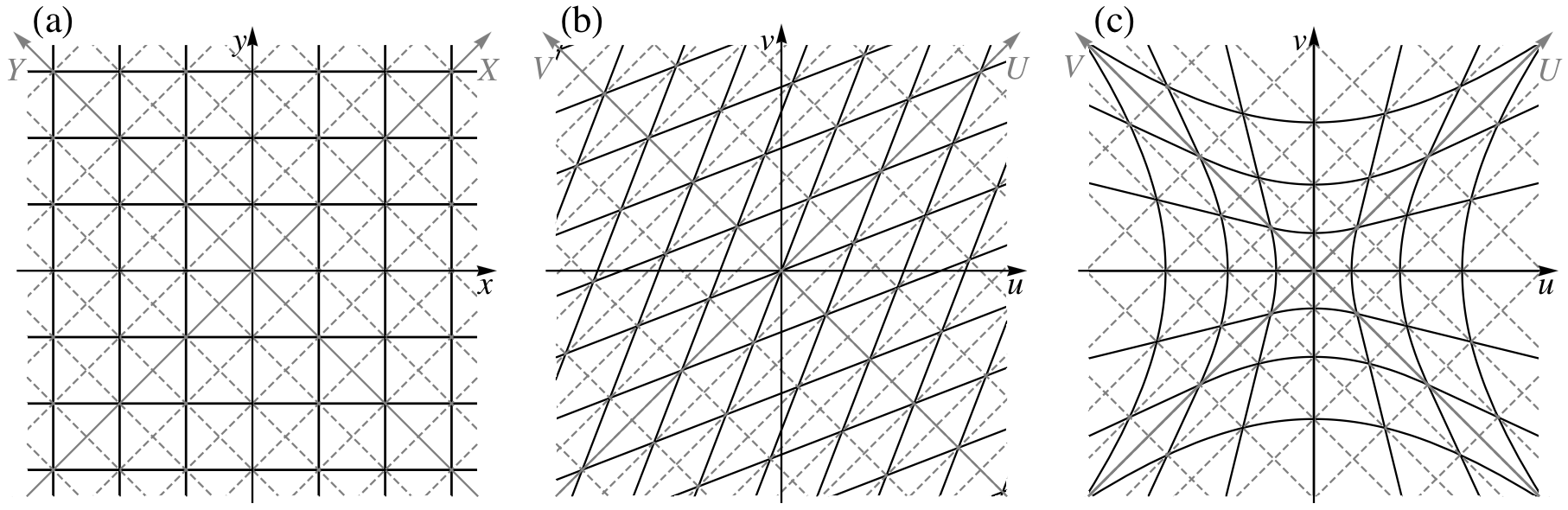}}
\caption{\small Coordinate curves of a linear (b) and a non-linear (c) mapping of the form $\alpha(X,Y) = (h(X), k(Y))$.  The bold curves are the images of the standard coordinate lines (constant $x$ or $y$), and the dotted lines are the images of the characteristic lines (constant $X$ or $Y$) illustrated in (a).  For (b), $h(t)=1.5 t, k(t)=\frac{2}{3}t$, and for (c), $h(t)=k(t)=\mbox{e}^t-1$. }
\label{fig:introduction}
\end{figure}

Such mappings send characteristic lines (constant $X$ or $Y$) to characteristic lines (constant $U$ or $V$).  A pair of horizontal and vertical lines (constant $x$ or $y$) are mapped to two curves whose tangent lines at the intersection point are reflections of each other about the characteristic lines, as shown in Fig.~\ref{fig:introduction}.

Each of the coordinate functions $U=h(X)$ and $V=k(Y)$ of the map $\alpha$ is itself a special solution of the wave equation.
Solutions of the wave equation and the maps that preserve it are immediately generalized by eliminating any requirement of smoothness and allowing the functions of the characteristic coordinates, $f_j$, $g_j$, $h$, $k$,  to be arbitrary.
This stands in stark contrast to the Laplace equation.
All of its distributional solutions are in fact analytic---they admit convergent power series.  The transformations that preserve them are conformal or anti-conformal, each component of which is a solution of the Laplace equation.

Whereas the Jacobian of a conformal map is a dilation-rotation, that of $(U,V)=(h(X),k(Y))$ is a variant of a hyperbolic dilation-rotation.  In $(x,y)$ and $(u,v)$ coordinates, this map and its cousin $(U,V)=(k(Y),h(X))$ are given by 
\begin{equation*}
  \textstyle (u,v) = \frac{1}{2}\big(\pm(h(x+y)-k(-x+y)), \ h(x+y)+k(-x+y)\big),
\end{equation*}
and both preserve the wave equation.  The Jacobians of these transformations have the form 
\begin{equation*}
  \renewcommand{\arraystretch}{1.0}
\left[
  \begin{array}{cc}
    \!u_x & u_y\! \\
    \!v_x & v_y\! \\
  \end{array}
\right] =
  \renewcommand{\arraystretch}{1.0}
\left[
  \begin{array}{cc}
    \!a & b\! \\
    \!b & a\! \\
  \end{array}
\right] \ \ \text{or} \ \ 
  \renewcommand{\arraystretch}{1.0}
\left[
  \begin{array}{rr}
    \!-a & -b\! \\
    \!b & a\! \\
  \end{array}
\right]\,,
\end{equation*}
respectively, each of which may preserve or reverse orientation.
In either case, we will call the transformation \textit{Lorentz conformal}, even when $h$ or $k$ are not differentiable.
See~\cite{Branson1987} for a discussion on Lorentz-conformal geometry in $n$-dimensional Minkowski space. 

Lorentz-conformal transformations include the linear Lorentz group, which preserves the Lorentz metric (see \cite{Naber}).  The reader will see in the pages that follow many explicit constructions of Lorentz-conformal transformations that realize certain geometric requirements.  One can envision a compendium of Lorentz-conformal mappings that complements the handbook of Moon and Spencer \cite{handbook} of Euclidean-conformal transformations.  
Triply orthogonal coordinate systems in dimension $2+1$ Lorentz space can be constructed from these planar mappings using operations of rotation, translation, and hyperbolic rotation
\cite{BP}.

\smallskip

Lorentz-conformal maps relate to the Lorentz inner product  analogously to the relationship of Euclidean conformal (and anti-conformal) maps to the Euclidean inner product.  
Consider a mapping $\alpha(x,y) = (u(x,y), v(x,y)) : \Omega \rightarrow \Omega'$, from a domain $\Omega \subseteq \mathbb{R}^{1,1}$ 
with the Lorentz inner product $\langle \ \  , \  \rangle = dx \otimes dx - dy \otimes dy$  
onto a region $\Omega' \subset \mathbb{R}^{1,1}$ with Lorentz inner product $\langle \ \  , \  \rangle' = du \otimes du - dv \otimes dv$.   
If $\alpha$ is differentiable and $v_x \neq 0 \neq u_y$ or $u_x \neq 0 \neq v_y$, then the following conditions are equivalent:

\begin{enumlist1}
  \item  $\alpha$ satisfies the following properties involving the Lorentz metric, in analogy to Euclidean conformal and anti-conformal mappings:
$$
\langle  \alpha_x, \alpha_y  \rangle' = u_x u_y - v_x v_y = 0,  \ \ \text{and} \ \ 
\langle  \alpha_x, \alpha_x  \rangle' = u_x^2 -v_x^2 = -u_y^2 +v_y^2 = -\langle \alpha_y, \alpha_y  \rangle'.
$$
Thus,
$\langle \ \  , \  \rangle' = du \otimes du - dv \otimes dv = \pm H^2 (dx \otimes dx - dy \otimes dy)$, where $H^2 =  |u_x^2 -v_x^2|$.
  \item  A Lorentzian analog of the Cauchy-Riemann equations holds:
$$
u_x = v_y, \; u_y = v_x
\qquad \text{or} \qquad
u_x = -v_y, \; u_y = -v_x\,.
$$
Therefore $u$ and $v$ satisfy the wave equation, $u_{xx}-u_{yy}=0$ and $v_{xx}-v_{yy}=0$.
\end{enumlist1}

If $\alpha$ is invertible, then (1) and (2) are equivalent to each of the following:

\begin{enumlist1}
\addtocounter{enumlist1}{2}
  \item  Both $\alpha$ and $\alpha^{-1}$ are Lorentz orthogonal; that is, for 
$c = (\det d\alpha)^{-2}$,
\begin{eqnarray*}
&& \langle \alpha_x, \alpha_y  \rangle' = u_x u_y - v_x v_y = 0  \\
&& \langle  \alpha_u^{-1}, \alpha_v^{-1}   \rangle = c(u_x v_x - u_y v_y) = 0\,.  
\end{eqnarray*}
  \item $\alpha^{-1}$ satisfies the following properties, in analogy to Euclidean conformal and anti-conformal mappings.   With $c = (\det d\alpha)^{-2}$,
\begin{eqnarray*}
&& \langle  \alpha_u^{-1}, \alpha_v^{-1}  \rangle = c(u_x v_x - u_y v_y)  = 0, \\
&& \langle \alpha_u^{-1}, \alpha_u^{-1}  \rangle = c(-v_x^2+v_y^2) = c(u_x^2-u_y^2) =  -\langle  \alpha_v^{-1}, \alpha_v^{-1}  \rangle\,.
\end{eqnarray*}
  \item  $\alpha$ preserves the wave equation: $f_{xx} - f_{yy} = 0$ if and only if $f_{uu}-f_{vv}=0$. 
\end{enumlist1}

Thus, Lorentz-conformal mappings are the Lorentzian analog of Euclidean conformal or anti-conformal maps (which together are sometimes called ``isothermal").  In contrast to the Euclidean case, 
where the mapping preserves or reverses orientation according to whether it is analytic or anti-analytic, 
mappings that satisfy either analog of the Cauchy-Riemann equations in Equivalence 2 may preserve or reverse orientation, and may preserve the $(1,-1)$ signature of the metric or reverse it to $(-1,1)$.  A single Lorentz-conformal mapping may even do both, on different parts of its domain, separated by lines of degeneracy where coordinate curves become tangent to each other.

Lorentz-conformal transformations may be used to simplify other hyperbolic equations.
An inhomogeneous version $f_{XY} = \nu(X) \mu(Y) g(f)$ of the nonlinear Klein-Gordon equation is brought into the standard form $f_{UV} = g(f)$ by the transformation $(U,V) = (h(X),k(Y))$ chosen so that $h'(X) = \nu(X)$ and $k'(Y) = \mu(Y)$.

\section{Lorentz-conformal transformation of lines and squares}\label{sec:geometric}

This section is devoted to discovering how simple geometric objects are transformed under invertible Lorentz-conformal transformations.
What kind of curve can a horizontal or vertical line be mapped to?  What are the possible images of a pair of intersecting horizontal and vertical lines?  What are all possible images of a square?

When one asks these questions for the inverse $\alpha^{-1}$ of a Lorentz-conformal transformation~$\alpha$, they become questions about the constant-$u$ and constant-$v$ contours of $\alpha$ itself.  We will take this point of view because it will be convenient for investigating transformations that are not one-to-one, but two-to-one or four-to-one, later on in Section~\ref{sec:symmetry}.

\smallskip

\subsection{Transformations of crossing lines}\label{sec:lines}

We investigate a basic question: Which pairs of crossing curves in the plane can be mapped to the coordinate axes under an invertible Lorentz-conformal transformation of the plane to itself, and what are all Lorentz-conformal transformations that do the job for a fixed pair of crossing curves?

\begin{figure}  
\centerline{
\scalebox{0.4}{\includegraphics{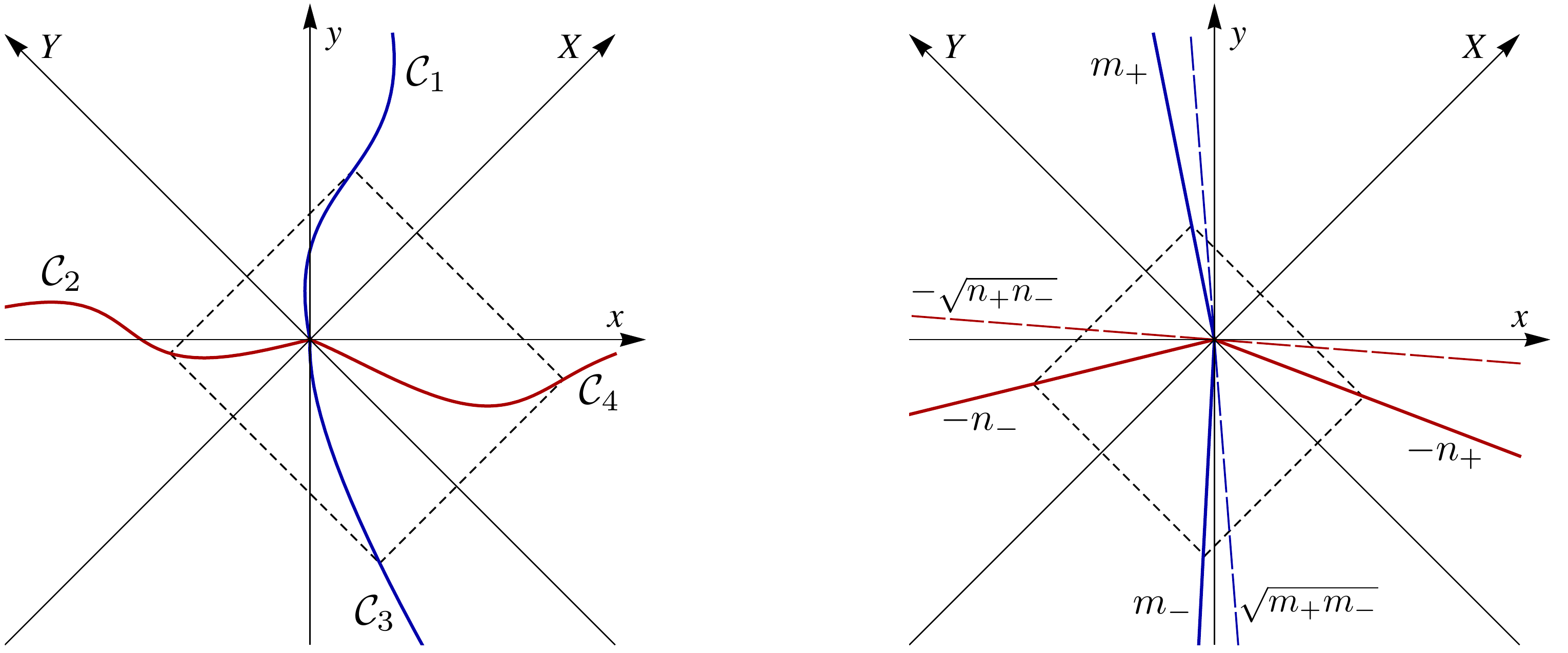}}
}
\caption{\small {\bfseries Left:}
If the curvilinear rays ${\cal C}_1$ and ${\cal C}_3$ together form the constant-$u$ contour of a Lorentz-conformal transformation from the $xy$ plane to the $uv$ plane and the rays ${\cal C}_2$ and ${\cal C}_4$ together form the constant-$v$ contour, then any three of the rays determines the fourth one by the rectangle rule (Rule~\ref{rule:rectangle1}).  
{\bfseries Right:}
The four tangent rays of the curves ${\cal C}_j$ at the origin also satisfy the rectangle rule.  The slope of the blue dashed line in characteristic coordinates $(X,Y)$ is the geometric average $\sqrt{m_+m_-}$ of the slope $m_+$ of ${\cal C}_1$ and $m_-$ of ${\cal C}_3$ at the origin, and the other dotted line is the analogous average of the the other two tangent lines.  These two averaged lines are reflections of each other about the $X$ and $Y$ axes.
}
\label{fig:rays}
\end{figure}

It turns out that any three of the four ``curvilinear rays" emanating from the intersection point (${\cal C}_j$ in Fig.~\ref{fig:rays}) determine the other one and that there is a whole functional degree of freedom in choosing a Lorentz-conformal transformation.  The description is pleasantly geometric, and the proofs involve neat functional equations involving monotonic functions.
Throughout Section~\ref{sec:geometric}, we deal with Lorentz-conformal mappings of the form $(U,V) = (h(X), k(Y))$ only, and not those of the form $(U,V) = (k(Y), h(X))$.

\smallskip
We start with the simple

\smallskip
\noindent {\bfseries Question:}  Which curves can be the constant-$u$ or constant-$v$ sets of an invertible Lorentz-conformal map?

\smallskip
Let $(h,k)$ be a pair of monotonic bijections from $\R$ to itself that define an invertible Lorentz-conformal map given in characteristic coordinates by $(U,V)=(h(X),k(Y))$.  The $u=u_0$ and $v=v_0$ contours are obtained from the change of coordinates (\ref{coordinates}),
\begin{eqnarray}
  && h(X) - k(Y) = 2u_0\,, \label{ucontour} \\
  && h(X) + k(Y) = 2v_0\,. \label{vcontour}
\end{eqnarray}
Each of these relations defines a monotonic bijective relation between $X$ and $Y$.  On the other hand, given a monotonic bijection $Y=g(X)$, one can set $h(X) := g(X)+2u$ and $k(Y):=Y$, and the resulting relation (\ref{ucontour}) is equivalent to $Y=g(X)$.  Such a relation is characterized geometrically by the property that it is continuous and intersects each characteristic line (constant $X$ or constant $Y$) exactly once.

Thus the answer to the question is

\begin{ruul}[Contours of invertible Lorentz-conformal maps]\label{rule:contours}
A set can be realized as the constant-$u$ or constant-$v$ set of an invertible Lorentz-conformal map $(U,V)=(h(X),k(Y))$, with $h$ and $k$ monotonic bijections from $\mathbb{R}$ to itself, if and only if it is a continuous curve $\cal C$ that intersects each characteristic line in exactly one point or, equivalently, it is a strictly monotonic 1-1 correspondence between the characteristic coordinates $X$ and~$Y$.

This condition, expressed in $(x,y)$ coordinates, is that $\cal C$ can be written either as 
$y=f(x)$ or $x=f(y)$, where the slope of any secant line of~$f$
cannot exceed $1$ in magnitude.
If $\cal C$ is differentiable, then $|f'|$ cannot exceed $1$ and cannot equal $1$ on any interval.
\end{ruul}

\begin{figure}  
\begin{center}
\scalebox{0.43}{\includegraphics{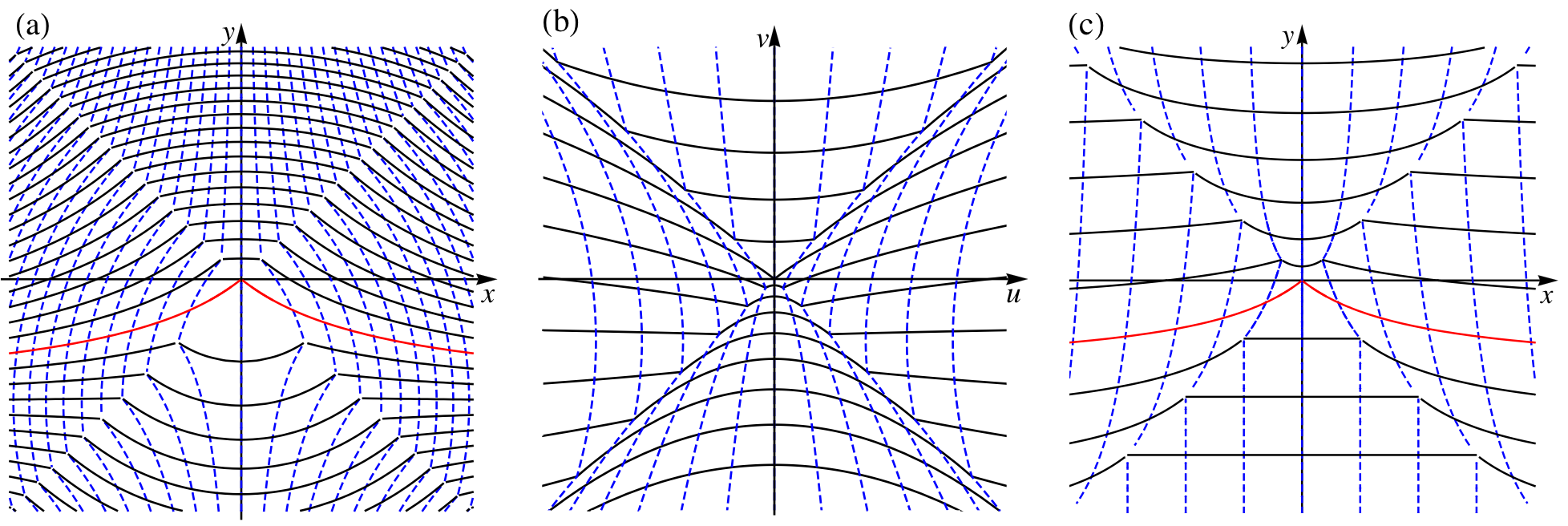}}
\end{center}
\setlength{\abovecaptionskip}{2pt plus 2pt minus 2pt}
\caption{\small Contour plots of Lorentz-conformal transformations in Examples~\ref{ex:ridge} and~\ref{ex:ridgeagain}.  The prescribed ridge curve (\ref{oceanridge}), highlighted in red in (a,c), is mapped to the $u$-axis.  (a) and (b) show contour plots of $\alpha (x,y) $ and $\alpha^{-1}(u,v)$ respectively, as given by (\ref{ridgeequations}).  (c) is a contour plot of $\alpha(x,y)$ given by~(\ref{equationsforridgeagain}).}
\label{fig:ridge}
\end{figure}

\begin{ex}[The ocean ridge]\label{ex:ridge}
Let us find a Lorentz-conformal map that has the ridge shape
$y = a/(1+|x|)$
as a contour (Fig.~\ref{fig:ridge}).  By Rule~\ref{rule:contours}, this is possible if and only if the constant~$a$ is no greater than $1$ in magnitude.  It is convenient to shift the relation so that the peak is at the origin:
\begin{equation} \label{oceanridge}
  y = \frac{-a|x|}{1+|x|}\,.
\end{equation}
In characteristic coordinates, this relation defines $Y$ as a decreasing function of $X$,
\begin{eqnarray}
  2(1+a)X + 2(1-a)Y + X^2 - Y^2 \,=\, 0 && \text{for $X\geq0$, $Y\leq0$}\,, \label{ridge1} \\
  2(1+a)Y + 2(1-a)X + Y^2 - X^2 \,=\, 0 && \text{for $X\leq0$, $Y\geq0$}\,, \label{ridge2}
\end{eqnarray}
which is realized as $h(X)+k(Y)=0$ with
\begin{equation} \label{ridgeequations}
  h(t) = k(t) =
\renewcommand{\arraystretch}{1.1}
\left\{
  \begin{array}{ll}
    2(1+a)t + t^2\,, & t\geq0\,, \\
    2(1-a)t - t^2\,, & t\leq0\,. 
  \end{array}
\right.
\end{equation}
Thus the ridge (\ref{oceanridge}) is the $v\!=\!0$ contour of the Lorentz-conformal map defined in characteristic coordinates by
$(U,V)=(h(X),k(Y))$.  Since $h\!=\!k$, the $u\!=\!0$ contour $h(X)-k(Y)=0$ is simply $X\!=\!Y$, or the $y$-axis.
\end{ex}

A characteristic feature of a linear Lorentz-conformal transformation (see Fig.~\ref{fig:introduction}~(b)) is that {\em any pair of coordinate lines $x = x_0$ and $y = y_0$ are mapped to lines that are reflections of one another about the characteristic lines}.  In other words, the slopes of the image lines are reciprocals of each other.   
This geometric feature carries over to the image of a pair of coordinate lines under a (nonlinear) differentiable Lorentz-conformal map, because the derivative of a Lorentz-conformal map at any point is a linear Lorentz-conformal transformation.  One simply takes the tangent lines at the point of intersection of the image curves, and these will have reciprocal slopes.  Applying this rule to a local inverse map shows that two intersecting contours of a Lorentz-conformal map have tangent lines that are reflections about the characteristic lines.

But what happens to this rule when a contour is not differentiable, such as in Example~\ref{ex:ridge}\,?

\smallskip
{\bfseries Question:} If a constant-$u$ contour and a constant-$v$ contour of a Lorentz-conformal map are left- and right-differentiable at their point of intersection, how are the four tangent lines at the point of intersection related?
\smallskip

The answer involves a certain average of the left and right tangent lines.
Define the {\em geometrically averaged tangent line} of a contour at a point $P$ to be the line through $P$ whose slope in characteristic coordinates $(X,Y)$ is the geometric average of the slopes 
of the contour's left and right tangent lines (preserving sign).

\begin{ruul}[\bfseries\slshape Nonsmooth crossing contours]\label{rule:nonsmoothcrossing}
Suppose two contours of an invertible Lorentz-conformal map are left- and right-differentiable at their crossing point $P$.
Then their four tangent lines at $P$ satisfy the following equivalent conditions, illustrated in Fig.~\ref{fig:rays}.
 \begin{enumlist}
 \item The geometrically averaged tangent lines of the two contours at $P$ are reflections of one another about either characteristic line through $P$.
 \item {\slshape (Rectangle rule)} If three of the tangent rays pass through vertices of a rectangle with sides parallel to the characteristic lines, then the fourth tangent ray must pass through the other vertex of the rectangle.
\end{enumlist}
\end{ruul}

The equivalence of (a) and (b) is seen with the aid of Fig~\ref{fig:rays}, which illustrates that condition (b) is equivalent to $m_-m_+=n_-n_+$, where $m_\pm$ and $n_\pm$ are the left and right derivatives in characteristic coordinates of intersecting contours.

Given that $h$ and $k$ are left- and right-differentiable at $X_0$ and $Y_0$, respectively, where $(X_0,Y_0)$ are the characteristic coordinates of $P$, one can prove this rule by differentiating the relations (\ref{ucontour}) and (\ref{vcontour}) implicitly.  Taking for simplicity $u_0=v_0=0$ as well as $h(0)=k(0)=0$, express these contours as the graphs of two monotonic functions $f$ and $g$:
\begin{equation*}
  \renewcommand{\arraystretch}{1.1}
\left.
  \begin{array}{llcl}
    u\equiv0: & h(X)-k(Y) = 0 & \iff & Y=f(X) \\
    v\equiv0: & h(X)+k(Y) = 0 & \iff & Y=g(X) 
  \end{array}
\right.
\end{equation*}
Differentiating each relation from the right and the left yields the system
\begin{equation*}
  \renewcommand{\arraystretch}{1.1}
\left[
  \begin{array}{cccc}
    0 & 1 & 0 & -f'(0+) \\
    1 & 0 & -f'(0-) & 0 \\
    0 & 1 & g'(0+) & 0 \\
    1 & 0 & 0 & g'(0-)
  \end{array}
\right]
\renewcommand{\arraystretch}{1.1}
\left[
  \begin{array}{c}
    h'(0-) \\ h'(0+) \\ k'(0-) \\ k'(0+)
  \end{array}
\right]  = 
\renewcommand{\arraystretch}{1.1}
\left[
  \begin{array}{c}
    0 \\ 0 \\ 0 \\ 0
  \end{array}
\right],
\end{equation*}
in which $\ell'(0\pm)$ denotes left and right derivatives at $0$.
The determinant of the matrix is $f'(0-)f'(0+) - g'(0-)g'(0+)=0$, which implies that the slopes of the geometrically averaged tangent lines of $f$ and $g$ at $0$ are equal in magnitude.

A simpler proof not involving $h$ and $k$ but only the crossing contours themselves will be given after Construction~\ref{const:realization1}, which provides functional relations between $(h,k)$ and the contours, thus eliminating the need for left- and right-differentiability of $h$ and $k$.

\begin{ex}[Smooth contours of a non-smooth map]\label{ex:nonsmooth}
This example shows that a constant-$u$ contour and a constant-$v$ contour can both be left- and right-differentiable at their crossing point $(X_0,Y_0)$ even if the functions $h$ and $k$ are not differentiable there.

Let $g:[1,2]\to[1,2]$ be any increasing bijection and define an increasing bijection $f:[0,\infty)\to[0,\infty)$ by $f(0)=0$ and, for positive integers $n$,
\begin{equation*}
  f(s) = 2^{-n} g(2^n s) \quad \text{ for }\; s\in[2^{-n},2^{-n+1}]\,.
\end{equation*}
Then define 
$h$ and $k$ by scaling the argument of $f$ by positive numbers $a_j$ as follows:
\begin{eqnarray}
  h(X) &=& U_0 + \renewcommand{\arraystretch}{1.1}
\left\{
  \begin{array}{rl}
    f(a_1(X-X_0)), & X\geq X_0\, , \\
    -f(a_3(X_0-X)), & X\leq X_0\,,
  \end{array}
\right. \\
  k(Y) &=& V_0 + \renewcommand{\arraystretch}{1.1}
\left\{
  \begin{array}{rl}
    f(a_2(Y-Y_0)), & Y\geq Y_0\,, \\
    -f(a_4(Y_0-Y)), & Y\leq Y_0\,.
  \end{array}
\right.
\end{eqnarray}
In characteristic coordinates $(X,Y)$ the contour $2u_0=h(X)-k(Y)$ is a straight line with slope $m_-=a_4/a_3$ for $X\leq X_0$ and a straight line with slope $m_+=a_2/a_1$ for $X\geq X_0$.  The contour $2v_0=h(X)+k(Y)$ has left and right slopes equal to $-n_-=-a_2/a_3$ and $-n_+=-a_4/a_1$.  These slopes satisfy $m_+m_-=n_+n_-$.

\end{ex}

It turns out that, if the two contours are left- and right-differentiable at their crossing point, then one can always find a Lorentz-conformal map that realizes these contours and for which $h$ and $k$ are differentiable from the right and left at the crossing.  We will return to this question after Construction~\ref{const:realization1}, which shows how to construct {\em all} Lorentz-conformal maps that realize a given pair of crossing contours.
The freedom in choosing the Lorentz-conformal map is balanced by a restriction on which contours can be realized.
This prompts the next question.

\smallskip
{\bfseries Question:} When can two intersecting curves be realized as the constant-$u$ and constant-$v$ contours of an invertible Lorentz-conformal map?

\smallskip
The answer is a generalization of the rectangle rule for the four tangent rays (Rule~\ref{rule:nonsmoothcrossing}b).  This is already clear in a simple case: If each half-contour emanating from the intersection point coincides with its tangent ray, then by the rectangle rule, any one of the rays would be determined by the other three.  In fact, this rule holds even when the rays are curvilinear, that is, one can choose three of the half-contours at will, and the other will be determined according to Fig.~\ref{fig:rays}.

\begin{ruul}[\bfseries\slshape Rectangle rule for crossing contours]\label{rule:rectangle1}
Four curvilinear rays\, ${\cal C}_j$ emanating from a point in $\R^2$ can be realized as the constant-$u$ and constant-$v$ contours of an invertible Lorentz-conformal map of the plane to itself if and only if both of the following conditions are satisfied.
  \begin{enumlist}
    \item There is one ray in each quadrant in characteristic coordinates $(X,Y)$, and each ray is a monotonic relation between $X$ and $Y$.
    \item If each of three rays contains one vertex of a given rectangle parallel to the characteristic lines, then the fourth ray contains the fourth vertex of the rectangle.  By this rule, any three rays determine the fourth (Fig.~\ref{fig:rays}).
  \end{enumlist}
\end{ruul}

This rule can be proved by simple geometric arguments.  Let $\alpha$ be a Lorentz-conformal transformation that takes the curvilinear rays ${\cal C}_j$ to the $u$ and $v$ coordinate rays, and let $\cal R$ be a rectangle with sides parallel to the characteristic lines and with three vertices on three of the rays ${\cal C}_j$.  Since $\alpha$ takes characteristic lines into characteristic lines, $\alpha({\cal R})$ is a rectangle with three vertices on three of the $uv$-coordinate rays.  Thus the fourth vertex of $\alpha({\cal R})$ is on the other coordinate ray.  This means that the fourth vertex of ${\cal R}$ lies the other of the four rays~${\cal C}_j$.

\smallskip
Another proof of the rectangle rule, involving functional equations, provides an explicit construction of the fourth curvilinear ray from the other three and an explicit construction of all invertible Lorentz-conformal maps $\alpha$ that realize a given pair of intersecting contours.  Knowing that both contours are completely determined by any three of the rays\, ${\cal C}_j$, one expects some freedom in the choice of $\alpha$.  Indeed, if $\lambda$ is an invertible Lorentz-conformal map that takes the lines $u=u_0$ and $v=v_0$ onto themselves, then the maps $\alpha$ and $\lambda\circ\alpha$ have the same $u=u_0$ and $v=v_0$ contours.  \label{lambda}
Supposing that $u_0\!=\!v_0\!=\!0$, such a transformation $\lambda$ maps both the standard coordinate axes and the characteristic coordinate axes onto themselves.  It possesses a full dihedral group $D_4$ of symmetries and requires that the associated functions $h$ and $k$ be equal to each other and odd, as we will see in Section~\ref{sec:symmetry}.  Thus, in the case that $h$ and $k$ are increasing, $\lambda$ is determined by a single increasing function from $[0,\infty)$ to itself.

\smallskip
Let ${\cal G}_{[0,\infty)}$ denote the group of increasing bijections from $[0,\infty)$ to $[0,\infty)$, and let ${\cal G}$ denote the group of increasing bijections from $\mathbb{R}$ to itself, under composition. 
Those invertible Lorentz-conformal transformations that take the $X$-axis to a line parallel to the $U$-axis and preserve the orientation of both of the characteristic axes form a group, which is isomorphic to ${\cal G}\times{\cal G}$.
A Lorentz-conformal map in this group is represented by $(h,k)\in{\cal G}\times{\cal G}$ through the relation $(U,V)=(h(X),k(Y))$:
\begin{equation*}
  (h,k)\in{\cal G}\times{\cal G}
  \quad\longleftrightarrow\quad
   (X,Y) \mapsto (U,V)=(h(X),k(Y))\,. 
\end{equation*}
For $(h,k)\in{\cal G}\times{\cal G}$, consider the $u$ and $v$ contours
\begin{eqnarray}
  && 2u_0 = h(X) - k(Y)\,, \quad ({\cal C}_1 \text{ and } {\cal C}_3)  \label{u0}\\
  && 2v_0 = h(X) + k(Y)\,. \quad ({\cal C}_2 \text{ and } {\cal C}_4)  \label{v0}
\end{eqnarray}
We may as well take the intersection point of these contours as well as the image thereof to be the origin: $(X_0,Y_0)=(0,0)$ and $(u_0,v_0)=(0,0)$.  The four curvilinear rays ${\cal C}_j$ emanating from the origin can be described by means of four increasing bijections of $[0,\infty)$.  This can be done in two ways: using the {\em positive and negative parts} of $h$ and $k$, 
\begin{equation}\label{hkplusminus}
\left.
  \begin{array}{l}
      h_+(s):=h(s) \\ h_-(s) := -h(-s) \\ k_+(s):=k(s) \\ k_-(s) := -k(-s)
  \end{array}
\right\},\,
s\geq0,
\end{equation}
or cyclically, by using four functions $g_j \in {\cal G}_{[0,\infty)}$ that express the curves ${\cal C}_j$ in characteristic coordinates:
\begin{equation}\label{Cg}
\renewcommand{\arraystretch}{1.1}
\left.
  \begin{array}{llcl}
    {\cal C}_1\,: & k_+(Y) = h_+(X) & \iff & \hspace{0.8em}Y = g_1(X) \,, \\
    {\cal C}_2\,: & h_-(-X) = k_+(Y) & \iff & -X = g_2(Y) \,, \\
    {\cal C}_3\,: & k_-(-Y) = h_-(-X) & \iff & -Y = g_3(-X) \,, \\
    {\cal C}_4\,: & h_+(X) = k_-(-Y) & \iff & \hspace{0.8em}X = g_4(-Y) \,.
  \end{array}
\right.
\end{equation}
The functions $g_j$ are determined by $h$ and $k$ through the map
\begin{equation}\label{ghkmap}
  (k_-,h_-,k_+,h_+) \;\mapsto\; (g_4,g_3,g_2,g_1) = (h_+^{-1}\circc k_-,\;k_-^{-1}\circc h_-,\;h_-^{-1}\circc k_+,\;k_+^{-1}\circc h_+)\,,
\end{equation}
which shows that they satisfy the cyclic condition
\begin{equation}\label{cyclic1}
  g_4\circ g_3\circ g_2\circ g_1 = \iota\,.
\end{equation}
In fact, this is exactly the solvability condition for equations~(\ref{ghkmap}).  To wit: given~$g_j$ ($j=1,2,3,4$), one can solve for three of the functions $(k_-,h_-,k_+,h_+)$ in terms of any one of them, say $k_-$, by using the last three coordinates of~(\ref{ghkmap}):
\begin{equation}\label{hkg1a}
\renewcommand{\arraystretch}{1.1}
\left.
  \begin{array}{lcl}
    k_- &=& p\,, \\
    h_- &=& p\circ g_3\,, \\
    k_+ &=& p\circ g_3\circ g_2\,, \\
    h_+ &=& p\circ g_3\circ g_2\circ g_1\,.    
  \end{array}
\right.
\end{equation}
The first coordinate of~(\ref{ghkmap}) is\, $k_-=h_+\circ g_4$, or
\begin{equation*}
  k_- = p\circ g_3\circ g_2\circ g_1\circ g_4\,
\end{equation*}
and is satisfied under condition (\ref{cyclic1}).
Therefore the functions (\ref{hkg1a}) determine a Lorentz-conformal map that realizes the curves ${\cal C}_j$ as constant $u$ and $v$ contours.

The cyclic condition (\ref{cyclic1}) on the functions $g_j$ is equivalent to the rectangle rule.
To wit: Consider a rectangle with sides parallel to the characteristic lines and having one vertex in each characteristic quadrant.  Let two opposite vertices have characteristic coordinates $(X_+,Y_+)$ and $(-X_-,-Y_-)$, with $X_\pm>0$, $Y_\pm>0$.
The passing of any ray ${\cal C}_j$ through the vertex in its quadrant can be expressed through the functions $g_j$ by
\begin{equation*}
\renewcommand{\arraystretch}{}
\left.
  \begin{array}{rcr}
     g_1(X_+)=Y_+ &\iff& (X_+,Y_+)\in{\cal C}_1\,, \\
    g_2(Y_+)=X_- &\iff&  (-X_-,Y_+)\in{\cal C}_2\,, \\
    g_3(X_-)=Y_- &\iff&  (-X_-,-Y_-)\in{\cal C}_3\,, \\
    g_4(Y_-)=X_+ &\iff&     (X_+,-Y_-)\in{\cal C}_4\,.
  \end{array}
\right.
\end{equation*}
Given any three of these, the other occurs if and only if $g_4\circ g_3\circ g_2\circ g_1 = \iota$.

\smallskip
The arbitrariness in the solution $(k_-,h_-,k_+,h_+)$ is in the free function $p$.  The set of solutions is the right coset of the diagonal subgroup
$\{ (p,p,p,p) : p \in{\cal G}_{[0,\infty)} \}$ of $\times_{j=1}^4{\cal G}_{[0,\infty)}$ that contains
$(\iota,\;g_3,\;g_3\circc g_2,\;g_3\circc g_2\circc g_1)$.
This coset can be expressed in terms of the functions $h$ and $k$ in $\cal G$ by extending $p$ to an odd increasing bijection $\ell$ from $\R$ to $\R$,
\begin{equation*}
  \ell(s) := \renewcommand{\arraystretch}{1}
\left\{
  \begin{array}{cl}
    p(s), & s\geq0\,,\\
    -p(-s), & s\leq0\,.
  \end{array}
\right.
\end{equation*}
One computes that the positive and negative parts of $\ell\circ h$ are $p\circ h_\pm$, and those of $\ell\circ k$ are $p\circ k_\pm$.
So $\ell\in{\cal G}$ realizes the Lorentz-conformal map $\lambda$ discussed above on page~\pageref{lambda}.

Equations (\ref{hkg1a}) provide a construction of $h$ and $k$ that utilizes the first three of the curves~${\cal C}_j$ and begins with choosing $k_-$ arbitrarily.  An unbiased way of viewing the construction is to rewrite (\ref{ghkmap}) to yield the identity
\begin{equation}\label{hkconstruction1}
  (h_+,k_-,h_-,k_+) \circ (g_4,g_3,g_2,g_1) = (k_-,h_-,k_+,h_+)\,,
\end{equation}
demonstrating how the $g_j$ cyclically permute the $(h,k)_\pm$.
For example, if one chooses $h_+$ arbitrarily, then $k_-$, $h_-$, and $k_+$ are determined successively by using $g_4$, $g_3$, and~$g_2$.

This construction and the freedom it allows are summarized in the following

\begin{const}[\bfseries\slshape Realizing crossing contours]\label{const:realization1}
\hspace{1em}
\begin{enumlist0}
\item
Let $\alpha$ be an invertible Lorentz-conformal map that realizes the curves ${\cal C}_j$, defined by $g_j$ in (\ref{Cg}), as $u$ and $v$ contours intersecting at the origin.  If $\alpha$ is represented by
$(U,V)=(h(X),k(Y))$, with $(h,k)\in{\cal G} \times {\cal G}$, the positive and negative parts $h_\pm$ and $k_\pm$ in ${\cal G}_0$ are determined by the cycle of equations
\begin{equation}\label{hkg1b}
  \renewcommand{\arraystretch}{1.1}
\left.
  \begin{array}{lcl}
  k_- &=& h_+\circ g_4\,, \\
  h_- &=& k_-\circ g_3\,, \\
  k_+ &=& h_-\circ g_2\,, \\
  h_+ &=& k_+\circ g_1\,.    
  \end{array}
\right.
\end{equation}
Exactly one of the functions $(k_-,h_-,k_+,h_+)$ can be chosen arbitrarily,
and the cyclic condition
\begin{equation}\label{cyclicRule1}
  g_4\circ g_3\circ g_2\circ g_1 = \iota\,,
\end{equation}
which is necessarily satisfied, guarantees the consistency of the four equations.
\item
Two invertible Lorentz-conformal maps $\alpha_1$ and $\alpha_2$ given in characteristic coordinates by \,$\{U=h_1(X),\,V=k_1(Y)\}$\, and \,$\{U=h_2(X),\,V=k_2(Y)\}$, respectively (with $(h_i,k_i)\in{\cal G} \times {\cal G}$), have the same $u\equiv0$ and $v\equiv0$ contours if and only if there is an increasing odd bijection $\ell$ of $\R$, with $\ell(0)=0$, such that
%
\begin{equation}\nonumber
  h_2=\ell \circ h_1 \quad\text{and}\quad k_2=\ell \circ k_1\,. 
\end{equation}
In other words, $\alpha_1$ and $\alpha_2$ are related by composition with a Lorentz-conformal transformation $\lambda$ defined by the pair $(h,k)=(\ell,\ell)$, with $\ell$ odd:
\begin{equation*}
  \alpha_2 = \lambda\circ\alpha_1\,.  
\end{equation*}
%

More generally, if $\alpha_1$ and $\alpha_2$ have the same $u\equiv u_0$ and $v\equiv v_0$ contours, one has $h_2(X) = \ell(h_1(X)-U_0)+U_0$ and $k_2(Y) = \ell(k_1(Y)-V_0)+V_0$, where $U_0=u_0+v_0$ and $V_0=-u_0+v_0$, or, more concisely,
\begin{equation*}
  \alpha_2 = \sigma\circ\lambda\circ\sigma^{-1}\,\circ\, \alpha_1\,,
\end{equation*}
in which $\sigma$ is a shift map in characteristic coordinates: $(U,V)\mapsto(U+U_0,V+V_0)$.
\end{enumlist0}
\end{const}

\smallskip
In light of the cyclic condition (\ref{cyclicRule1}) on the contours, one can give a simpler proof of Rule~\ref{rule:nonsmoothcrossing}.
If their right and left derivatives 
exist at their crossing point, the $g_j$ are right and left differentiable at $0$.  In fact, if three of them are differentiable, then all of them are.  Differentiating the cyclic condition yields
\begin{equation*}
  g_4'(0)g_3'(0)g_2'(0)g_1'(0) = 1.
\end{equation*}
If the slopes of the curvilinear rays ${\cal C}_j$ in characteristic coordinates are as given in Fig.~\ref{fig:rays}, one has
$m_+=g'_1(0)$, $m_-=g'_3(0)$, $n_+=1/g'_4(0)$ and $n_-=1/g'_2(0)$ and thus obtains
\begin{equation*}
 m_-m_+ = n_-n_+\,,
\end{equation*}
which shows that the geometrically averaged tangent lines of the two contours are reflections of each other about the characteristic axes.

The first proof of Rule~\ref{rule:nonsmoothcrossing}, given after the statement of the rule, was valid whenever $h$ and $k$ were left-and right-differentiable at the crossing point.  But in fact, Construction~\ref{const:realization1} provides enough flexibility in choosing the transformation that realizes the two crossing contours so that $h$ and $k$ can be chosen to be left- and right-differentiable.  To see this, let us revisit Example~\ref{ex:nonsmooth}.

\begin{ex}[Example~\ref{ex:nonsmooth} revisited]\label{ex:nonsmoothagain}
The two intersecting contours contours are left- and right-differentiable at their crossing point although the $h$ and $k$ for the transformation are not.  The function $f$ that imparts the non-differentiability to $h$ and $k$ plays the role of the function $p$ in 
(\ref{hkg1a}).

If $p$ is chosen to be smooth, then one obtains new functions $h$ and $k$ that are left- and right-differentiable and generate a Lorentz-conformal map that realizes the same intersecting contours.
In the case that $p=\iota$, one obtains (for $U_0=V_0=X_0=Y_0=0$) functions $(h,k)_\pm$ and $g_j$ that are linear.
\begin{equation} \nonumber 
\renewcommand{\arraystretch}{1.1}
\left.
  \begin{array}{lcl}
    h_+(t) = a_1 t &\quad& g_1(t) = (a_1/a_2)t \\
    h_-(t) = a_3 t &\quad& g_2(t) = (a_2/a_3)t \\
    k_+(t) = a_2 t &\quad& g_3(t) = (a_3/a_4)t \\
    k_-(t) = a_4 t &\quad& g_4(t) = (a_4/a_1)t \\
  \end{array}
\right.
\end{equation}
\end{ex}

Let us see what Construction~\ref{const:realization1} has to say about the ocean ridge Example~\ref{ex:ridge}, where there was a natural polynomial choice of $h$ and $k$.  From the point of view of the cyclic equations (\ref{hkg1b}), or, equivalently, the mapping (\ref{ghkmap}), these choices of $h$ and $k$ are not so immediate.

\begin{ex}[Example~\ref{ex:ridge} revisited]\label{ex:ridgeagain}
The curves ${\cal C}_2$ and ${\cal C}_4$ are prescribed by the ridge shape.  Its reflection symmetry about the $y$-axis is equivalent to $g_2 = g_4^{-1}$.  In fact, whenever this holds, ${\cal C}_1$ and ${\cal C}_3$ can be taken to be the positive and negative rays of the $y$-axis, which corresponds to $g_1=g_3=\iota$ because this satisfies $g_4\circc g_3\circc g_2\circc g_1=\iota$.  It follows from Construction~\ref{const:realization1} (\ref{hkg1b}) that $h=k$.  

One finds the relations $X=g_4(-Y)$ and $-X=g_2(Y)$ in Example~\ref{ex:ridge} by setting $h(X)+k(Y)=0$ (the contour $v\equiv0$), or equivalently from (\ref{ridge1},\ref{ridge2}):
\begin{eqnarray}
  g_2(t) &=& \left[ (t+(1+a))^2 - 4a \right]^{1/2} - (1-a)\,, \\
  g_4(t) = g_2^{-1}(t) &=& \left[ (t+(1-a))^2 + 4a \right]^{1/2} - (1+a)\,.
\end{eqnarray}
Assuming $g_1=g_3=\iota$, if one sets $k_-=\iota$, then Construction~\ref{const:realization1} gives $k_-=h_-=\iota$ and $k_+=h_+=g_2$, or
\begin{equation}\label{equationsforridgeagain}
  h(t) = k(t) =
\renewcommand{\arraystretch}{1.1}
\left\{
  \begin{array}{ll}
    \left[ (t+(1+a))^2 - 4a \right]^{1/2} - (1-a)\,, & t\geq0\,, \\
    t\,, & t\leq0\,.
  \end{array}
\right.
\end{equation}
Another possible choice of $h$ and $k$ comes from taking $h_+=\iota$,
\begin{equation}\nonumber
  h(t) = k(t) =
\renewcommand{\arraystretch}{1.1}
\left\{
  \begin{array}{ll}
    t\,, & t\geq0\,, \\
   \left[ (t+(1-a))^2 + 4a \right]^{1/2} - (1+a)\,, & t\leq0\,.
  \end{array}
\right.
\end{equation}
The second option is obtained from the first by composition with the odd function
\begin{equation*}
  \ell(t) = {\rm{sgn}}(t)\,g_4(|t|)\,.
\end{equation*}
\end{ex}

\subsection{Transformations of a square}\label{sec:square}

\smallskip
{\bfseries Question:} What shapes are mapped by an invertible Lorentz-conformal transformation to a square parallel to the coordinate axes?
\smallskip

\begin{figure}  
\centerline{\scalebox{0.4}{\includegraphics{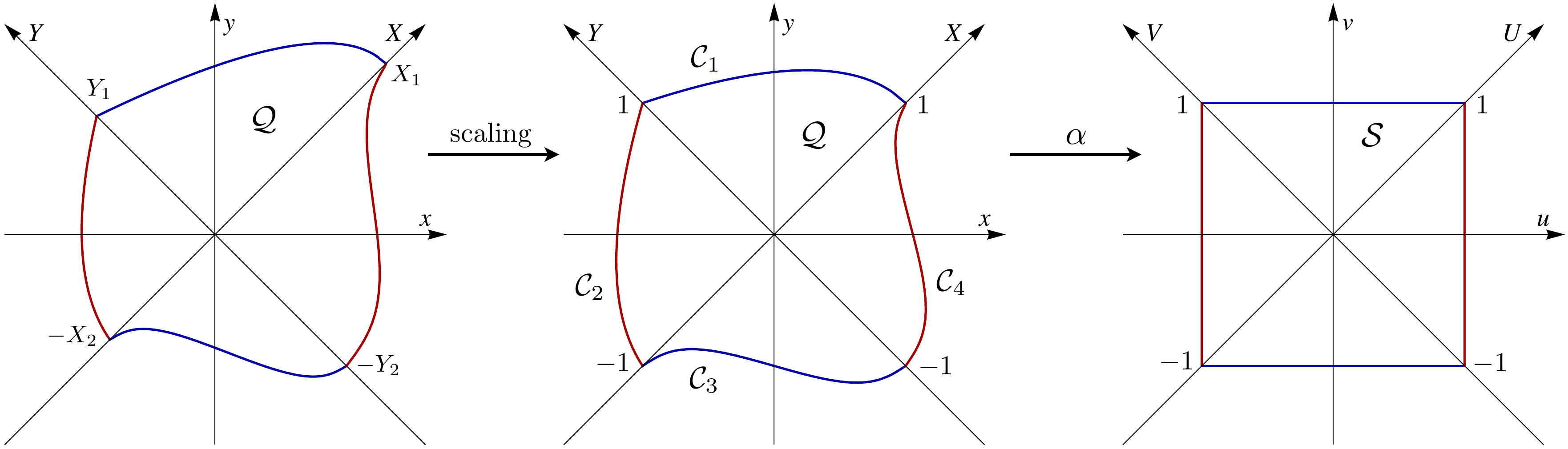}}}
\caption{\small A curvilinear quadrilateral with vertices on the characteristic axes (at $X_1$, $-X_2$, $Y_1$, $-Y_2$) can be transformed Lorentz-conformally to a quadrilateral whose vertices are those of a square by scaling the positive and negative characteristic ($X$ and $Y$) axes by appropriate positive numbers.  A Lorentz-conformal transformation $\alpha$ will take this new curvilinear quadrilateral $\cal Q$ to the standard square ${\cal S}$ as long as the four sides ${\cal C}_j$ are related by the rectangle rule (Rule~\ref{rule:rectangle2}; see Fig.~\ref{fig:quadrilateralB}).}
\label{fig:quadrilateralA}
\end{figure}

\begin{figure}  
\centerline{\scalebox{0.4}{\includegraphics{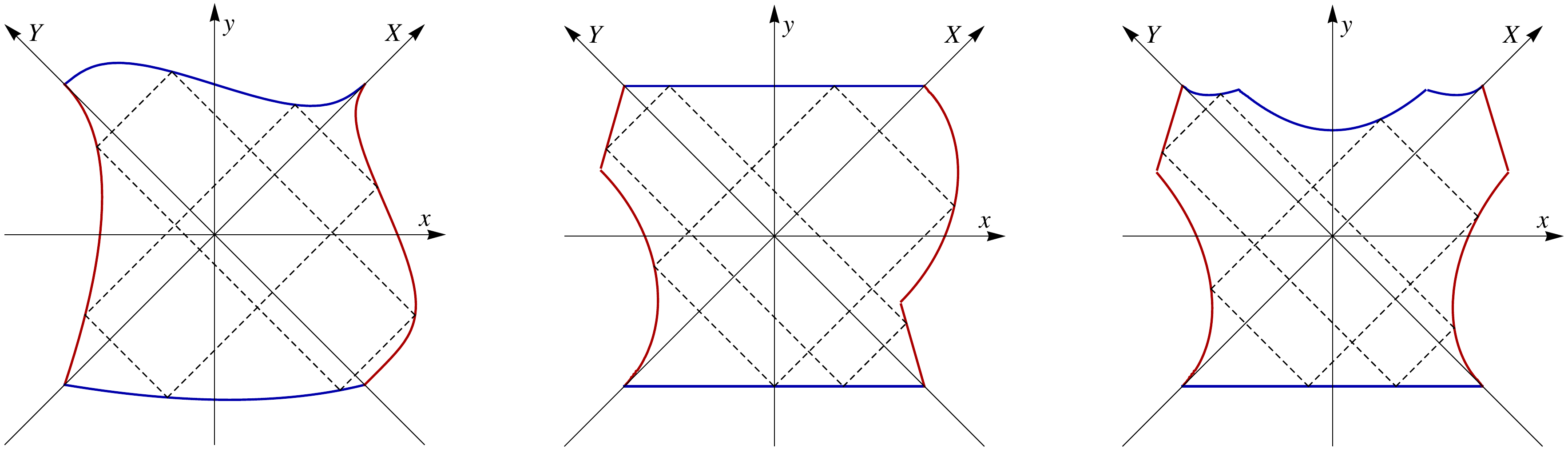}}}
\caption{\small Each of these three curvilinear quadrilaterals $\cal Q$ can be mapped 
Lorentz-conformally to the standard square (Fig.~\ref{fig:quadrilateralA}, right) because they obey the rectangle rule (Rule~\ref{rule:rectangle2}).  Any rectangle with three vertices on different sides of the quadrilateral must have its fourth vertex on the other side of the quadrilateral.
{\bfseries Center:} If the top and bottom of $\cal Q$ are straight lines, then any lateral side is obtained from the other by reflection across the $x$-axis and translation (Construction~\ref{rule:flattopbottom}).  {\bfseries Right:} If the bottom of $\cal Q$ is flat and the lateral sides are reflections of each other about the $y$-axis, then the top is symmetric about the $y$-axis (Construction~\ref{rule:symmetricsides}).}
\label{fig:quadrilateralB}
\end{figure}

Denote the transformation by $\alpha$.
The diagonals of the square are characteristic lines in the $(u,v)$ plane, which are mapped by $\alpha^{-1}$ back to crossing characteristic lines in the $(x,y)$ plane.  The sides of the square are mapped by $\alpha^{-1}$ to constant-$u$ and constant-$v$ curve segments that intersect each other at four points on these crossing characteristic lines (Fig.~\ref{fig:quadrilateralA}).
The curve segments form a ``curvilinear quadrilateral"; each side remains within its own quadrant formed by the two characteristic lines.
Figures~\ref{fig:quadrilateralA} and \ref{fig:quadrilateralB} shows some curvilinear quadrilaterals $\cal Q$ that can be Lorentz-conformally transformed into a square $\cal S$.  Similarly to the problem of intersecting contours, it turns out that only three of the sides can be determined arbitrarily.  The fourth is determined from the other three by a rectangle rule.

This problem is natural is if one is trying to solve the wave equation in a bounded domain by mapping it to a square, where solutions can be obtained by separation of variables.  Thinking of $y$ as time, a curvilinear quadrilateral $\cal Q$ may represent a spatial interval with moving endpoints.  If both variables are spatial, one may be considering two-dimensional internal gravity waves
(see \cite{Mercier2008,MartinSmith2011}, for example).

\begin{ruul}[\bfseries\slshape Rectangle rule for transformations of a square]\label{rule:rectangle2}
A curvilinear quadrilateral can be mapped Lorentz-conformally to a square with sides parallel to the coordinate axes under and only under the following conditions.
\begin{enumlist}
\item The vertices lie on a pair of crossing characteristic lines.
\item Each side intersects any characteristic line in at most one point; equivalently, the sides are monotonic in characteristic coordinates.
\item If each of three sides contains one vertex of a given rectangle parallel to the characteristic lines, then the fourth side contains the fourth vertex of the rectangle.  By this rule, any three sides determine the fourth, as illustrated in Fig.~\ref{fig:quadrilateralB}.
\end{enumlist}
\end{ruul}
\noindent

Part (c) of this rule has a simple geometric proof that is similar to that of Rule~\ref{rule:rectangle1}.  Let $\cal Q$ be a curvilinear quadrilateral such that $\alpha({\cal Q})={\cal S}$, where ${\cal S}$ is the right square centered at the origin and with sides of length 2 in the $uv$ plane shown in Fig.~\ref{fig:quadrilateralB}(right).
Let $\cal R$ be a rectangle with sides parallel to the characteristic lines and three vertices on three sides of~$\cal Q$.  Since $\alpha$ maps characteristic lines into characteristic lines, $\alpha({\cal R})$ is a rectangle, and three of its vertices lie on three sides of the square~$\cal S$.  Thus the fourth vertex of $\alpha({\cal R})$ lies on the fourth side of~$\cal S$, so the fourth vertex of~$\cal R$ lies on the fourth side of~$\cal Q$.

There is a nice interpretation of this rule if one thinks of $(x,y)$ as denoting space and time and the lateral sides of $\cal Q$ as moving boundaries.  Choose any point on the bottom side of $\cal Q$, and send out a signal in each direction with unit celerity.  One signal hits the left side of $\cal Q$ and bounces back, and the other hits the right side and bounces back.  
Where these two signals meet again determines a point on the top of $\cal Q$.  If this is done for each point on the bottom, all the meeting points will trace out the top side of $\cal Q$.

\smallskip

A more elaborate proof provides an explicit construction of the fourth side from the other three and reveals the structure of all Lorentz-conformal transformations that realize a mapping from a given $\cal Q$ to $\cal S$.
It makes use of the group ${\cal G}_{[-1,1]}$ of increasing bijections of $[-1,1]$ to itself and the group ${\cal G}_{[0,1]}$ of increasing bijections from $[0,1]$ to itself, with the multiplication being composition; the identity is denoted by~$\iota$.  There are two involutions of ${\cal G}_{[0,1]}$ that commute with each other: the anti-isomorphism $g\mapsto g^{-1}$, which reflects the graph of $g$ about the $45^\circ$ line; and the isomorphism that rotates the graph of $g$ by $180^\circ$ about the center point of $[0,1]^2$, which is denoted by $g\mapsto \tilde g$\,:
\begin{equation*}
  \tilde g(s) = 1-g(1-s)\,.
\end{equation*}
%
%
%

The proof is similar to that of Rule~\ref{rule:rectangle1}.
First, by translating a curvilinear quadrilateral $\cal Q$ so that its vertices lie on the characteristic axes and then applying a simple Lorentz-conformal transformation that scales the four characteristic rays through the origin, $\cal Q$ may be transformed so that its vertices coincide with those of the unit square $\cal S$, or $X_1=X_2=Y_1=Y_2=1$ in Fig~\ref{fig:quadrilateralA}.  We will thus assume, for the remainder of this section, that $\cal Q$ has this property.

A Lorentz-conformal transformation $\alpha$ given by $(U,V) = (h(X),k(Y))$, with $h,k \in {\cal G}_{[-1,1]}$, maps the diamond $(X,Y)\in[-1,1]^2$ containing $\cal Q$ bijectively to the diamond $(U,V)\in[-1,1]^2$ containing $\cal S$.  It can be extended to an invertible transformation of $\R^2$ by extending the domains of $h$ and $k$ to all of $\R$.
Suppose that $\alpha$ takes ${\cal Q}$ to $\cal S$.
In terms of $h$ and $k$, the curves ${\cal C}_j$ are described~by
\begin{equation*}
  \renewcommand{\arraystretch}{1.2}
\left.
  \begin{array}{ll}
    {\cal C}_1\,: & h(X)+k(Y)=1, \\
    {\cal C}_2\,: & h(X)-k(Y)=-1, \\
    {\cal C}_3\,: & h(X)+k(Y)=-1, \\
    {\cal C}_4\,: & h(X)-k(Y)=1,
  \end{array}
\right.
\quad X, Y \in [-1,1].
\end{equation*}
By evaluating at the intersection points of pairs of curves, one finds that $h(-1)=k(-1)=-1$, $h(0)=k(0)=0$, and $h(1)=k(1)=1$.  
Thus the positive and negative parts of $h$ and $k$, namely the functions $h_\pm$ and $k_\pm$ defined by (\ref{hkplusminus}), are increasing bijections of $[0,1]$ ({\itshape i.e.,} they are in ${\cal G}_{[0,1]}$).  The curves ${\cal C}_j$ are described in two ways,
\begin{equation}\label{hkpm}
  \renewcommand{\arraystretch}{1.2}
\left.
  \begin{array}{lrcr}
    {\cal C}_1\,: & h_+(X)+k_+(Y)=1 &\iff& Y = g_1(1-X)\,, \\
    {\cal C}_2\,: & h_-(-X)+k_+(Y)=1 &\iff& -X = g_2(1-Y)\,, \\
    {\cal C}_3\,: & h_-(-X)+k_-(-Y)=1 &\iff& -Y = g_3(1+X)\,, \\
    {\cal C}_4\,: & h_+(X)+k_-(-Y)=1 &\iff& X = g_4(1+Y)\,,
  \end{array}
\right.
\quad X, Y \in [0,1],
\end{equation}
in which the functions $g_j$ are in ${\cal G}_{[0,1]}$.

The $g_j$ are related to $h_\pm$ and $k_\pm$ by the mapping
\begin{equation*}
  (k_-,h_-,k_+,h_+) \;\mapsto\; (g_4,g_3,g_2,g_1) = (h_+^{-1}\circc\tilde k_-,\;k_-^{-1}\circc\tilde h_-,\;h_-^{-1}\circc\tilde k_+,\;k_+^{-1}\circc\tilde h_+)\,,
\end{equation*}
and they satisfy the cyclic condition
\begin{equation}\label{cyclic2}
  g_4\circ\tilde g_3\circ g_2\circ\tilde g_1 = \iota\,,
\end{equation}
which is in fact the solvability condition for $(k_-,h_-,k_+,h_+)$ in terms of the $g_j$.
The solution has one degree of freedom $p \in{\cal G}_{[0,1]}$,
\begin{equation}\label{hkg2a}
\renewcommand{\arraystretch}{1.1}
\left.
  \begin{array}{lcl}
  k_- &=& \tilde p\,,\\
  h_- &=& p\circ\tilde g_3\,,\\
  k_+ &=& \tilde p\circ g_3\circ\tilde g_2\,,\\
  h_+ &=& p\circ\tilde g_3\circ g_2\circ\tilde g_1\,.    
  \end{array}
\right.  
\end{equation}
More compactly, one can solve for $(h,k)_\pm$ recursively with the equations
\begin{equation}\label{hkconstruction2}
  (h_+,k_-,h_-,k_+) \circ (g_4,g_3,g_2,g_1) = (\tilde k_-,\tilde h_-,\tilde k_+,\tilde h_+)\,,
\end{equation}
by choosing one of the four arbitrarily and solving cyclically for the other three.  Condition~(\ref{cyclic2}) will guarantee that the four equations are consistent.

Condition (\ref{cyclic2}) is equivalent to the rectangle rule.
To wit:
The conditions that the vertices $(X_+,Y_+)$, $(-X_-,Y_+)$, $(-X_-,-Y_-)$, and $(X_+,-Y_-)$ lie on the curves ${\cal C}_1$, ${\cal C}_2$, ${\cal C}_3$, ${\cal C}_4$ are expressed by
\begin{eqnarray*}
  (X_+,Y_+)\in {\cal C}_1 &\Longleftrightarrow& g_1(1-X_+)=Y_+\,, \\
  (-X_-,Y_+)\in {\cal C}_2 &\Longleftrightarrow& \tilde g_2(Y_+)=1+X_- \,, \\
  (-X_-,-Y_-)\in {\cal C}_3 &\Longleftrightarrow& g_3(1+X_-)=-Y_-\,, \\
  (X_+,-Y_-)\in {\cal C}_4 &\Longleftrightarrow& \tilde g_4(-Y_-)=1-X_+\,,
\end{eqnarray*}
respectively.
The first three of these conditions, for example, implies $g_3\circ\tilde g_2\circ g_1(1-X_+)=-Y_-$.  Thus the fourth condition holding for all $X_+\in[0,1]$ is equivalent to the condition~(\ref{cyclic2}).

\begin{const}[\bfseries \slshape Realizing transformations to a square]\label{const:realization2}
\hspace{1em}
\begin{enumlist0}
\item
Let $\alpha$ be an invertible Lorentz-conformal transformation that maps the curvilinear quadrilateral defined by the curves ${\cal C}_j$, $j=1,\dots,4$, as in the right-hand-side of (\ref{hkpm}), to the unit square about the origin (Fig.~\ref{fig:quadrilateralA}).  If $\alpha$ is represented by $(U,V)=(h(X),k(Y))$, with $h$ and $k$ in ${\cal G}_{[-1,1]}$,
the functions $h$ and $k$ are determined by their positive and negative parts $h_\pm$ and $k_\pm$ in ${\cal G}_{[0,1]}$ through the cycle of equations
\begin{equation}\label{hkg2b}
\renewcommand{\arraystretch}{1.1}
\left.
  \begin{array}{lcl}
  \tilde k_- &=& h_+\circ g_4\,, \\
  \tilde h_- &=& k_-\circ g_3\,, \\
  \tilde k_+ &=& h_-\circ g_2\,, \\
  \tilde h_+ &=& k_+\circ g_1\,.    
  \end{array}
\right.  
\end{equation}
Exactly one of the functions $(k_-,h_-,k_+,h_+)$ can be chosen arbitrarily, and the cyclic condition
\begin{equation}\label{cyclicRule2}
  g_4\circ\tilde g_3\circ g_2\circ\tilde g_1 = \iota\,,
\end{equation}
which is necessarily satisfied, guarantees the consistency of the four equations.
\item
Two invertible Lorentz-conformal maps $\alpha_1$ and $\alpha_2$ 
given in characteristic coordinates by \,$\{U=h_1(X),\,V=k_1(Y)\}$\, and \,$\{U=h_2(X),\,V=k_2(Y)\}$, respectively,
map the same curvilinear quadrilateral 
${\cal Q}$ to the unit square ${\cal S}$ (Fig.\,\ref{fig:quadrilateralA}) if and only if there is an increasing odd bijection $\ell$ of $[-1,1]$ with $\ell(0)=0$, such that, on $[-1,1]$,
\begin{equation}\nonumber
  h_2 = \ell\circ h_1
  \qquad
  \text{and}
  \qquad
  k_2 = \tilde\ell\circ h_2\,,
\end{equation}
in which $\tilde\ell$ is defined by $\tilde\ell(s)=1-\ell(1-s)$ for $s\in[0,1]$ and is extended to an odd function on $[-1,1]$.
In other words, $\alpha_1$ and $\alpha_2$ are related, on the domain $[-1,1]^2$ in $XY$ coordinates, by composition with a Lorentz-conformal transformation $\lambda$ defined by the pair $(h,k)=(\ell,\tilde\ell)$:
\begin{equation*}
  \alpha_2 = \lambda\circ\alpha_1\,.
\end{equation*}
\end{enumlist0}
\end{const}

\smallskip
Suppose one desired to map a quadrilateral ${\cal Q}$ with flat top and bottom, as in Fig.~\ref{fig:quadrilateralB} (center), to the standard square ${\cal S}$.  This is natural if one is thinking of $y$ as representing time and the region ${\cal Q}$ as representing, for example, a string with endpoints moving at a celerity less than $1$.
The sides ${\cal C}_1$ and ${\cal C}_3$ are represented by the identity function $g_1\!=\!\iota$ and $g_3\!=\!\iota$.  By the cyclic relation $g_4\circc\tilde g_3\circc g_2\circc\tilde g_1 = \iota$, one has $g_4\circc g_2=\iota$.
This condition is equivalent to the geometric condition that ${\cal C}_2$ is obtained from ${\cal C}_4$ by reflection about the $x$-axis and translation.

Let us set $g_2=g$ and $g_4=g^{-1}$ for some $g\in{\cal G}_{[0,1]}$.
Relations (\ref{hkg2a}) or (\ref{hkg2b}) show that all possible transformations $\alpha$ 
(represented by $(h,k) \in {\cal G}_{[-1,1]} \times {\cal G}_{[-1,1]}$) that map ${\cal Q}$ to ${\cal S}$ are given by
\begin{equation*}
  (k_-,h_-,k_+,h_+) = (\tilde p,\; p,\; \tilde p\circc\tilde g,\; p\circc g)\,,
\end{equation*}
in which $p\in{\cal G}_{[0,1]}$ is arbitrary.  One obtains the following rule.

\smallskip
\begin{const}[Flat top and bottom]\label{rule:flattopbottom}
A curvilinear quadrilateral ${\cal Q}$ with flat top and bottom is mapped by a Lorentz-conformal transformation to the unit square ${\cal S}$ if and only if the two lateral sides of ${\cal Q}$ are related by reflection about the $x$-axis and translation.  After translation and scaling the characteristic rays so that the vertices of $\cal Q$ coincide with those of ${\cal S}$ (Fig.\,\ref{fig:quadrilateralA}), the Lorentz-conformal transformations of the form $(U,V) = (h(X), k(Y))$
that realize such a mapping are represented exactly by those pairs $(h,k)$ whose positive and negative parts in ${\cal G}_{[0,1]}$ are related by $h_\pm = \tilde k_\pm$.
  
The functions $g_j$ describing the sides ${\cal C}_j$ of ${\cal Q}$ are related to $(h,k)$ by
\begin{equation*}
  (k_-,h_-,k_+,h_+) = (\tilde p,\; p,\; \tilde q,\; q)
  \;\xrightarrow{g \;=\; p^{-1}\circ\, q}\;
  (g_1,\;g_2,\;g_3,\;g_4) = (\iota,\;g,\;\iota,\;g^{-1})
\end{equation*}
The right-hand side $(\iota,\;g,\;\iota,\;g^{-1})$ characterizes all ${\cal Q}$, with vertices coinciding with those of ${\cal S}$,
whose top and bottom are flat and that can be mapped to ${\cal S}$ by a Lorentz-conformal transformation, and the left-hand side characterizes the Lorentz-conformal transformations that map such ${\cal Q}$ to ${\cal S}$.
\end{const}

\begin{figure}[h]  
\centerline{\includegraphics[width=5.9in]{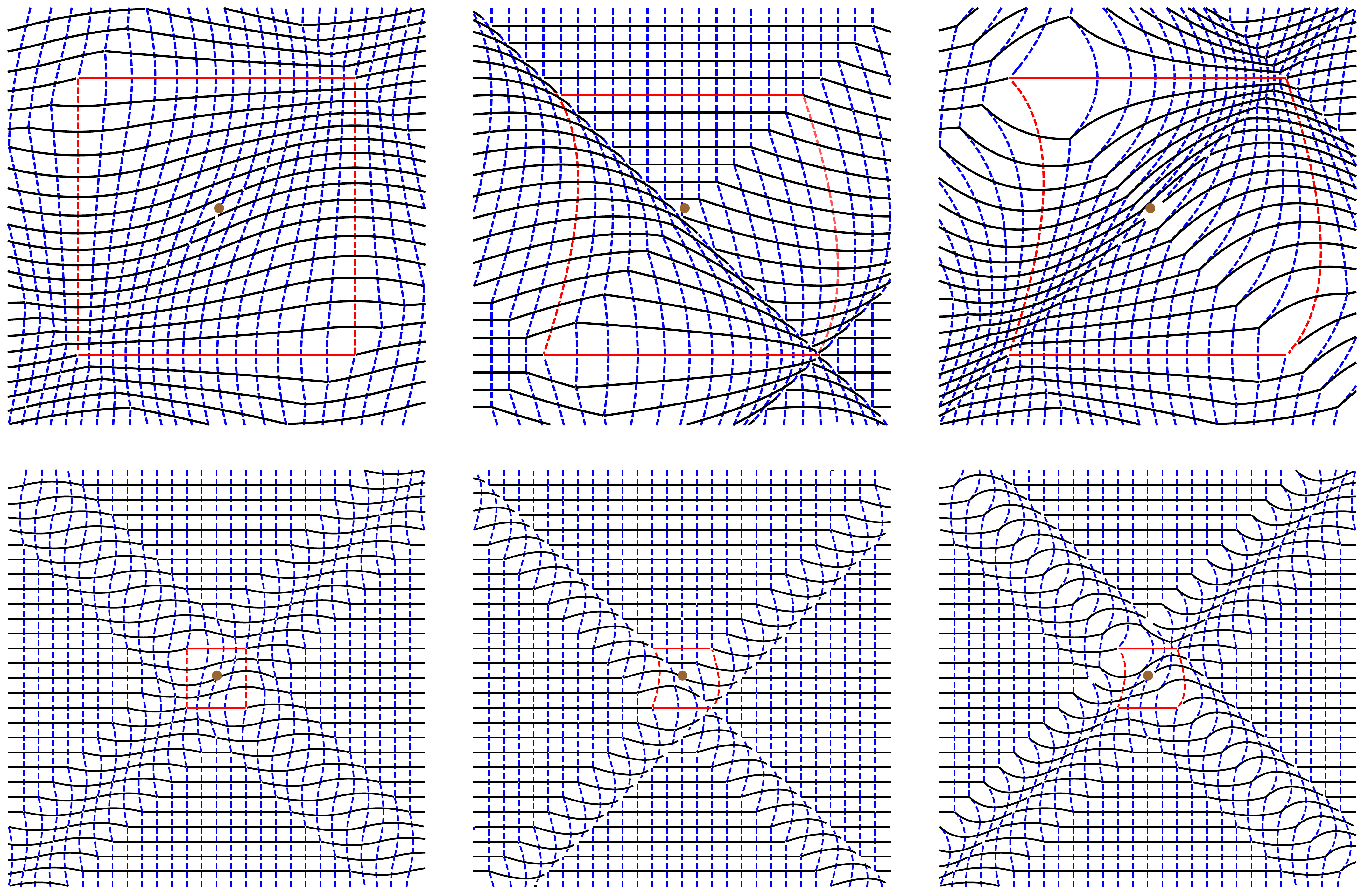}}
\label{fig:flat1}
\caption{\small  Contour plots in the $(u,v)$-plane of Lorentz-conformal transformations illustrating Const.~\ref{rule:flattopbottom}. For $|t|\leq 1$, $g(t)$ and $p(t)$ are as follows: \textbf{Left Column:} $g(t)=t,p(t)=(\mbox{e}^t-1)/(e-1)$, \textbf{Center Column:} $g(t)=t^2,  \  p(t)=\sqrt{t},$ \textbf{Right Column:} $g(t)=t^2, \ p(t)=(\mbox{e}^t-1)/(e-1)$.  For $|t| \geq 1$, $g(t) = p(t) = t$. \textbf{Top Row:} $(u,v) \in [-0.75,0.75] \times [-0.75,0.75]$, \textbf{Bottom Row:} $(u,v) \in [-3.5,3.5] \times [-3.5,3.5]$.  The origin is marked by a (brown) dot. }
\end{figure}

Now let us demand that the bottom of ${\cal Q}$ be flat and that the lateral sides be symmetric under reflection about the $y$-axis, as in Fig.~\ref{fig:quadrilateralB} (right).  This means that $g_2=g$ and $g_4=\tilde g^{-1}$ for some $g\in{\cal G}_{[0,1]}$.  The cyclic relation reduces to $g_1 = \tilde g^{-1}\circc g$, and therefore
\begin{equation*}
  \tilde g_1^{-1} = g_1\,,
\end{equation*}
which is equivalent to the rotational symmetry of the top side ${\cal C}_1$ about the $y$-axis.
The Lorentz-conformal transformations that map ${\cal Q}$ to ${\cal S}$ are given by
\begin{equation*}
  (k_-,h_-,k_+,h_+) = (\tilde p,\; p,\; \tilde p\circc\tilde g,\; p\circc\tilde g)\,.
\end{equation*}

\begin{const}[Left-right symmetry]\label{rule:symmetricsides}
  If a curvilinear quadrilateral ${\cal Q}$  whose bottom is flat and whose lateral sides are reflections of one another about the $y$-axis is mapped by a Lorentz-conformal transformation $\alpha$ to the unit square ${\cal S}$, then the top side is symmetric about the $y$-axis.  If $\alpha$ is given by $(U,V)=(h(X),k(Y))$, with $h$ and $k$ in ${\cal G}_{[-1,1]}$, then the functions $g_j$ in (\ref{hkpm}) describing the sides ${\cal C}_j$ of ${\cal Q}$ are related to $(h,k)$ by
\begin{equation*}
  (k_-,h_-,k_+,h_+) = (\tilde p,\; p,\; \tilde p\circc p^{-1}\circc q,\; q)
  \;\xrightarrow{g \;=\; \tilde p^{-1}\circ\, \tilde q}\;
  (g_1,\;g_2,\;g_3,\;g_4) = (\tilde g^{-1}\circc g,\;g,\;\iota,\;\tilde g^{-1})
\end{equation*}
The right-hand side $(\tilde g^{-1}\circc g,\;g,\;\iota,\;\tilde g^{-1})$ characterizes all ${\cal Q}$ 
(with vertices coinciding with those of ${\cal S}$)
whose bottom is flat and whose lateral sides are reflections of one another about the $y$-axis, and the left-hand side characterizes all Lorentz-conformal transformations $(U,V) = (h(X), k(Y))$ that map such ${\cal Q}$ to ${\cal S}$.
\end{const}

\begin{figure}[h]  
\centerline{\includegraphics[width=5.9in]{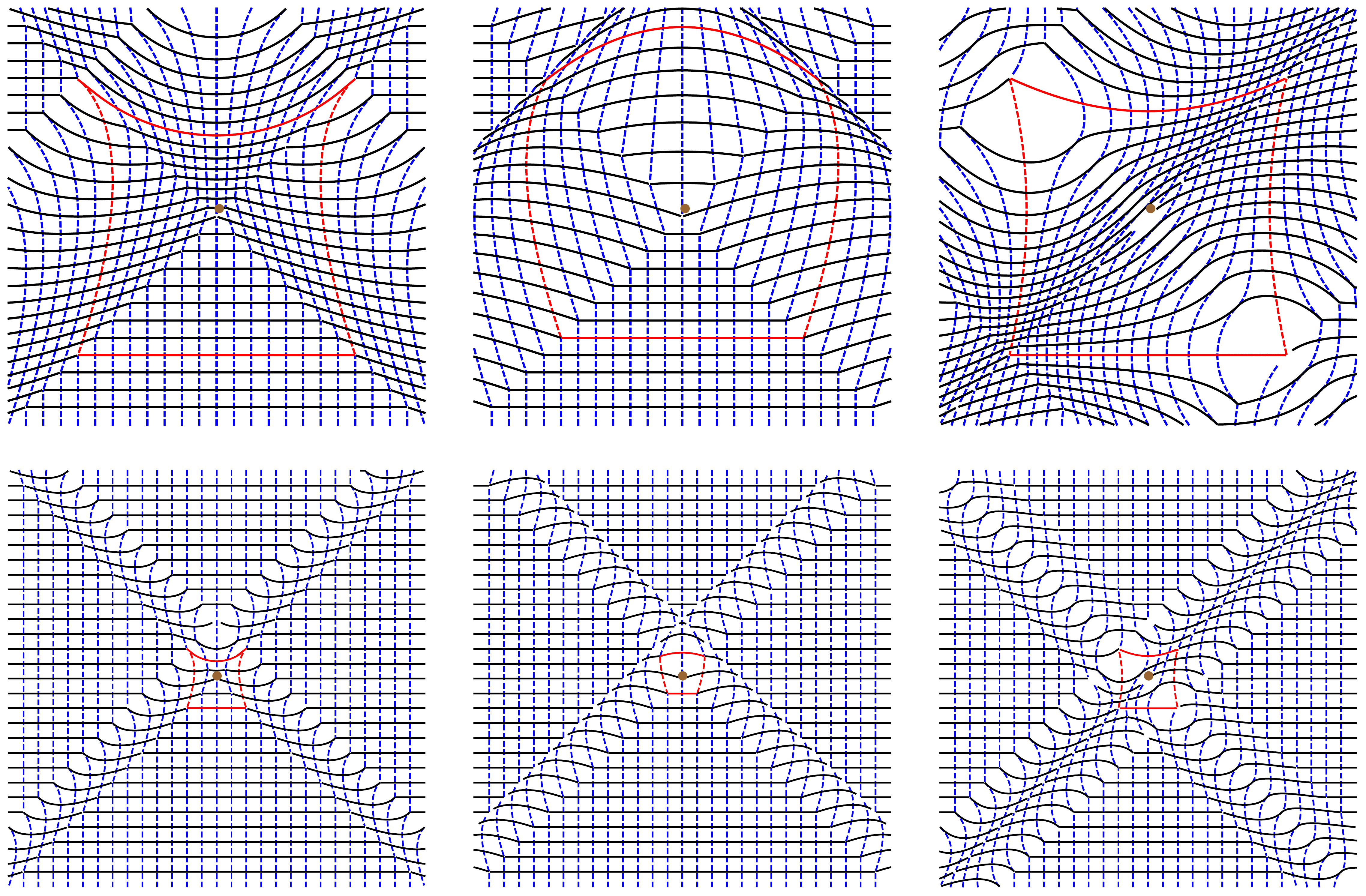}}
\label{fig:flat2}
\caption{\small Contour plots in the $(u,v)$-plane of Lorentz-conformal transformations illustrating Const.~\ref{rule:symmetricsides}. For $|t| \leq 1$, the functions $g(t)$ and $p(t)$ are as follows: \textbf{Left Column:} $g(t)=t^2, p(t)=t$, \textbf{Center Column:} $g(t)=\sqrt{t}, p(t)=t$, \textbf{Right Column:} $g(t)=(\mbox{e}^t-1)/(e-1), p(t)=t^2$.  For $|t| \geq 1$, $g(t) = p(t) = t$. \textbf{Top Row:} $(u,v) \in [-0.75,0.75] \times [-0.75,0.75]$, \textbf{Bottom Row:} $(u,v) \in [-3.5,3.5] \times [-3.5,3.5]$.  The origin is marked by a (brown) dot.  }  
\end{figure}

\section{Symmetry and unfolding}\label{sec:symmetry}

In how many ways can one color the pattern in Fig.~\ref{Section3IntroductionFigure}(a) or (b) or (c) with two colors to create the constant-$u$ and constant-$v$ contours of a Lorentz-conformal transformation?  For the first two patterns, there are four such colorings, and each corresponds to a Lorentz-conformal map with different symmetries.  In this section, we ``unfold" many-to-one Lorentz-conformal mappings into invertible ones in such a way that the contour plots of the two transformations are different colorings of one another.   
We  characterize distinguished classes of Lorentz-conformal mappings in terms of symmetries of their contour plots under the dihedral group $D_4$.  Unfoldings of many-to-one mappings are reflected in a change of the symmetry group.

\centerline{\includegraphics[width=5.9in]{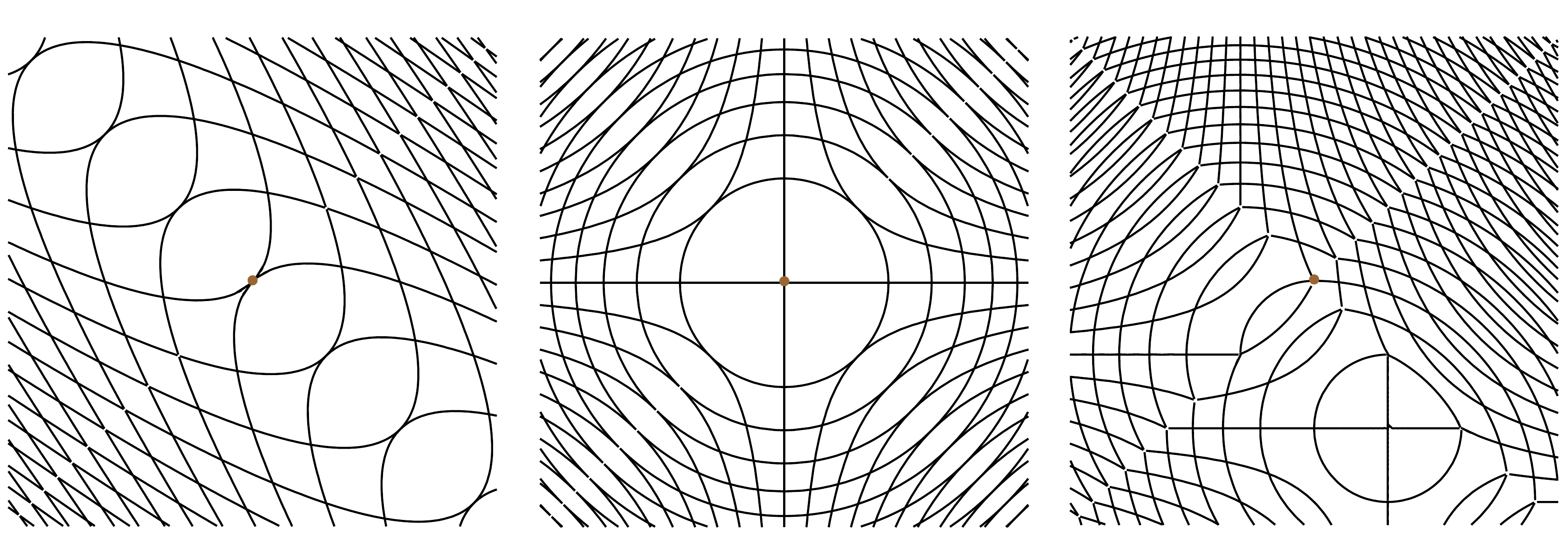}}

\begin{figure}[h]
\setlength{\abovecaptionskip}{3pt}
\caption{\small Contour plots of Lorentz-conformal mappings, without distinguishing the two families of coordinate curves by color.  In how many ways can each pattern be colored with two colors to create the contour plot of a Lorentz-conformal mapping?  The (brown) dot marks the origin.}
\label{Section3IntroductionFigure}
\end{figure}

\subsection{Unfolding a 4:1 mapping and cropping out degeneracies}\label{unfoldingsection}

For a non-invertible function $h : \mathbb{R} \rightarrow \mathbb{R}$, the transformation $\beta$ given in characteristic coordinates by $(U,V) = (h(X), h(Y))$ is Lorentz-conformal, but not invertible.  If $h : \mathbb{R} \rightarrow [0, \infty)$ is even and its restriction to $[0, \infty)$ is bijective, then $\beta$ is 4:1 on the complement of the characteristic axes and maps each quadrant cut out by these axes onto the first quadrant $U >0$, $V>0$ (Fig.~\ref{maptowedge}).  Writing $(u,v) = \beta(X,Y)$, with
$$u = \frac{1}{2}(h(X) - h(Y)) \ \ \mbox{and} \ \ v = \frac{1}{2}(h(X) + h(Y)),$$
we see that for $h$ even, reflection across $X = 0$ or $Y = 0$ 
folds every constant-$u$ and constant-$v$ contour onto itself.  

\smallskip
{\bfseries Question:} How can the 4:1 mapping $\beta$, with $h$ even, be ``unfolded" to create a 1:1 Lorentz-conformal transformation?
\smallskip

By restricting $h$ to $[0, \infty)$ and extending it antisymmetrically to $\mathbb{R}$, Construction~\ref{unfolding} creates an invertible Lorentz-conformal transformation $\alpha$ that has the same contour plot as $\beta$ in the quadrant
$X>0$, $Y>0$ and whose full contour plot is a different coloring of the contour plot of $\beta$.
Under this unfolding, constant-$u$ contours of $\beta$ are cut and re-joined with constant-$v$ contours of $\beta$ at each characteristic axis to create the contour curves of the unfolded  
 transformation $\alpha$.  This is shown in Fig.~\ref{unfoldingquadratic}(a,b), with $h(t) = t^2$.  The unfolded transformation, and its inverse, have the property that 
reflection about a characteristic axis interchanges constant-$u$ and constant-$v$ curves, rather than preserving them, as does the 4:1 mapping $\beta$.  
We say that an $n$:1 Lorentz-conformal mapping $\alpha$ is an {\it unfolding} of an $m$:1 mapping $\beta$ with $n<m$ if the contour plots of $\alpha$ and $\beta$ coincide on some fundamental domain where $\beta$ is invertible,
and look the same when all contours are depicted in the same color.  

\begin{const}[An unfolding]\label{unfolding}
Let $p : \mathbb{R} \rightarrow [0, \infty)$  be an even function whose restriction $p_+ : [0, \infty) \rightarrow [0, \infty)$ is bijective, so that $\beta$ given by $(U,V) = (p(X), p(Y))$ is 4:1 on the complement of the characteristic axes.    Let
$h : \mathbb{R} \rightarrow \mathbb{R}$ be the odd extension of $p_+$ to $\mathbb{R}$:
\begin{equation} \label{hoddextension} 
h(t) = 
\left\{ \hspace{-5pt} \begin{array}{rr} p(t), & t \geq 0, \\ 
-p(t),  &  t < 0. \end{array} \right.
\end{equation}
The  Lorentz-conformal transformation $\alpha$ given by $(U,V) = (h(X),h(Y))$ is an unfolding of~$\beta$.  Restricting $p$ 
to $(-\infty,0]$ and extending it antisymmetrically to $\mathbb{R}$ creates a variation of this unfolding, $(U,V)=-(h(X),h(Y))$.
\end{const}
\noindent 

\begin{ex}[Unfolding a quadratic mapping]\label{quadratic}
Consider the 4:1 map $\beta$ given by $(U,V) = (X^2, Y^2)$, or $(u,v) = (2xy, x^2 + y^2)$.  For every point $(u_0,v_0)$ in the interior of the image of $\beta$ (the quadrant $U>0, V>0$), the contour $u=u_0$ (hyperbola) intersects the contour $v=v_0$ (circle) in four points (Fig.~\ref{maptowedge}). 
The contours of the unfolded transformation $\alpha$ from
Construction~\ref{unfolding} 
are obtained from the contours of $\beta$ by cutting and re-joining circles with hyperbolas at every point on a characteristic axis (Fig.~\ref{unfoldingquadratic}(a,b)).  The unfolded transformation is also an example of Construction 2: the circle (a curvilinear quadrilateral) is mapped to a square.
\end{ex} 

\centerline{\includegraphics[width=5.0in]{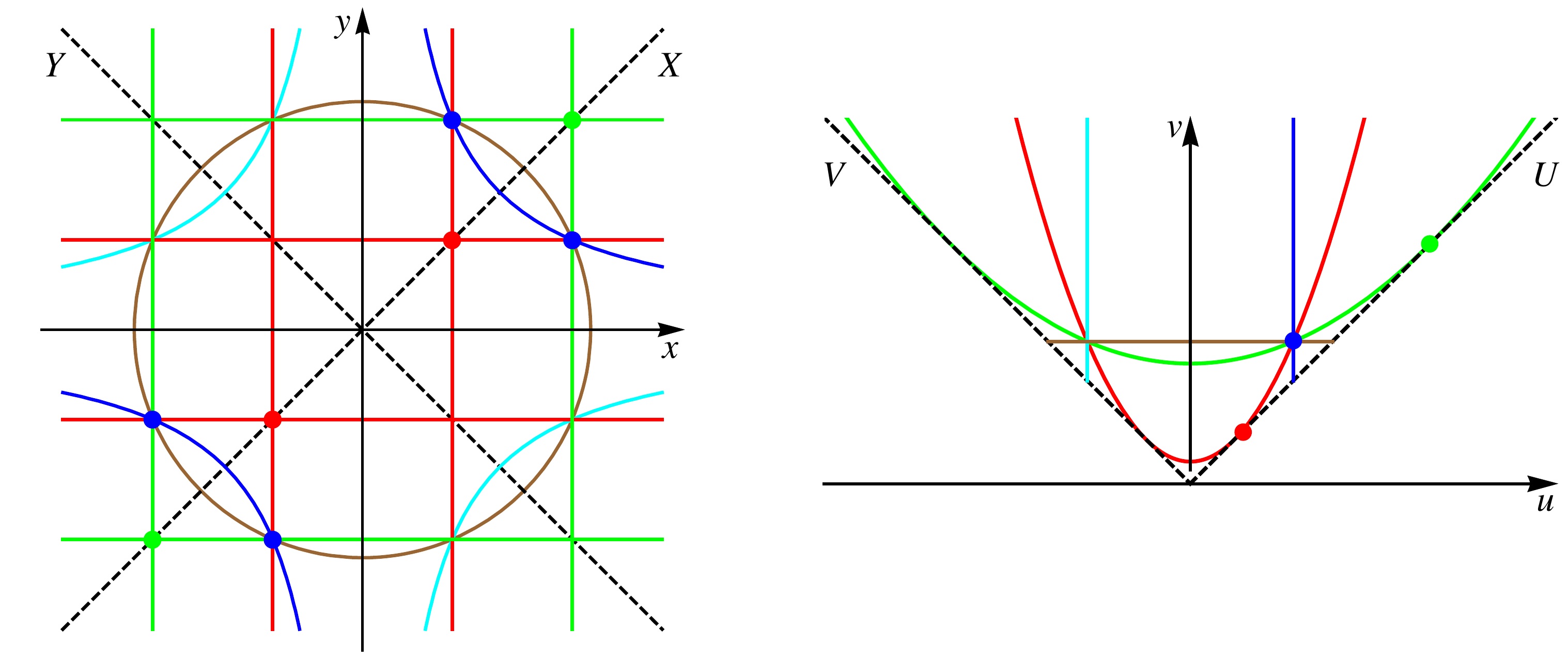}}

\begin{figure}[h]
\setlength{\abovecaptionskip}{2pt}
\caption{\small The 4:1 mapping $(u,v) = (2xy, x^2 + y^2)$ in Example~\ref{quadratic}, which maps each quadrant of the $(X,Y)$-plane onto the first quadrant of the $(U,V)$-plane. 
}
\label{maptowedge}
\end{figure}

A property of unfolding is that the resulting transformation, when it is invertible, is a change of coordinates in that any two coordinate curves intersect each other exactly once, rather than multiple times (see Fig.~\ref{unfoldingquadratic}(a,b)).  Folding occurs in Euclidean-conformal mappings as well, as in some of the coordinate plots shown in \cite{handbook}.  For example, the analytic function $w = z^2$ gives the 2:1 conformal mapping $(u,v) = (x^2-y^2, 2xy)$, which is the Euclidean analog of the 4:1 Lorentz-conformal mapping
$(u,v) = (2xy, x^2 + y^2)$ of Example~\ref{quadratic}.  

\centerline{\includegraphics[width=6in]{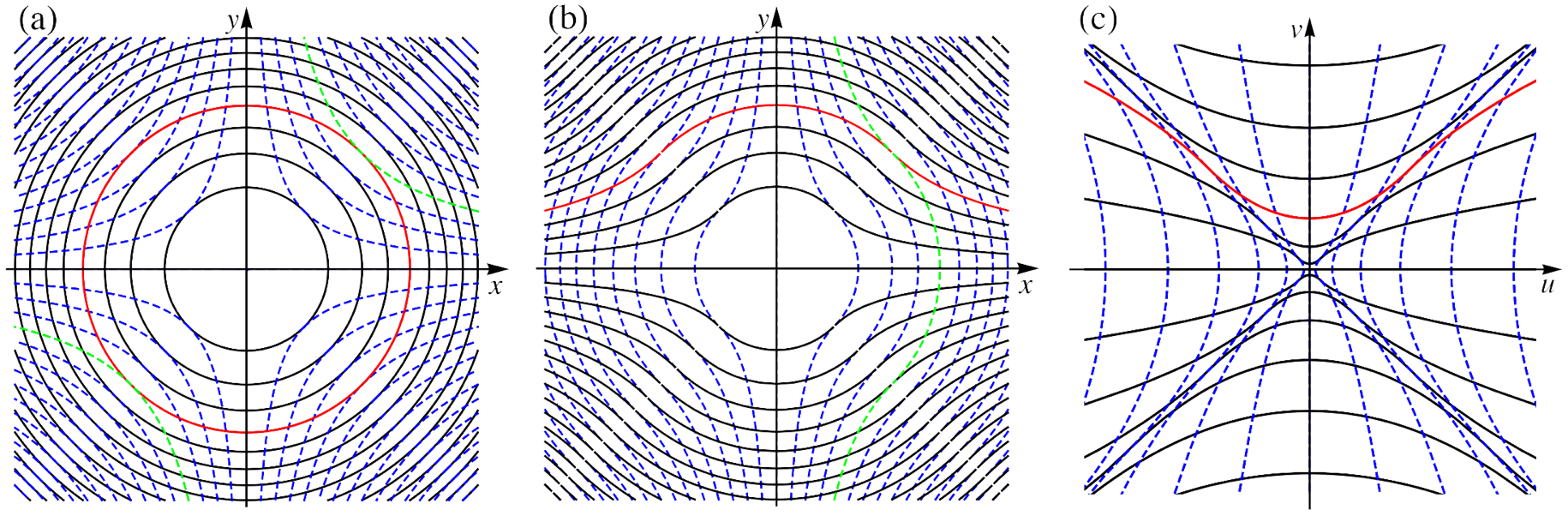}}

\begin{figure}[h]
\caption{\small (a) Contour plot of the 4:1 Lorentz-conformal mapping $\beta$ given by $(u,v) = (2xy, x^2 + y^2)$. (b) Contour plot of the bijective mapping $\alpha$, which unfolds $\beta$ via Construction~\ref{unfolding}.  (c)  $x$ and $y$ contours of $\alpha^{-1}$.  In (a) and (b), red marks one $v$-contour  and green marks one $u$-contour.  In (c), red marks one $x$-contour.}
\label{unfoldingquadratic}
\end{figure}

Lorentz-conformal mappings (including invertible ones) may have degeneracies in which the differential is not surjective, as in Fig.~\ref{unfoldingquadratic}, where $u$ and $v$ contours intersect tangentially at the characteristic axes.  In a true global coordinate system, contours intersect transversely everywhere, as in the exponential and logarithmic plots in Fig.~\ref{explogplots} (b,c). 

\bigskip

\centerline{\includegraphics[width=6in]{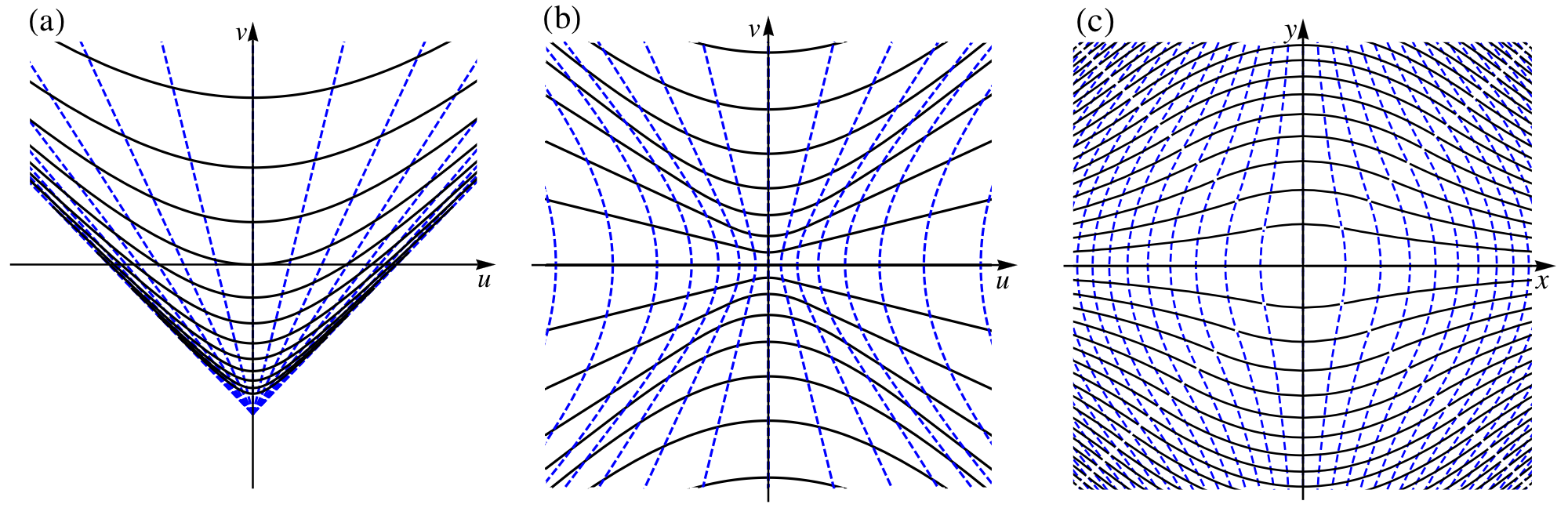}}
\begin{figure}[h]
\setlength{\abovecaptionskip}{2pt}
\caption{\small
Constant $x$ and $y$ contours of (a): $\beta$: $(U,V) = (p(X), p(Y))$ with $p(t) = e^t -1$, (b): $\alpha$ from Construction~\ref{unfolding} with $p_+(t) = e^t -1$.  (c): Constant $u$ and $v$ contours of $\alpha^{-1}$. }
\label{explogplots}
\end{figure}

Suppose each of $h$ and $k$ is differentiable, except possibly at points where the tangent line to its graph is vertical.  We ask the following

\smallskip
{\bfseries Question:} What is the locus of degeneracy where constant $u$ and $v$ contours of $\alpha$: $(X,Y) = (h(X),k(Y))$ are tangent to each other?
\smallskip

Since Lorentz-orthogonal vectors are reflections of each other about characteristic lines, the $u$ and $v$ contours of $\alpha$ through a point $(X,Y)= (P,Q)$ are tangent if and only if their common tangent line at $(P,Q)$ is in a characteristic direction.  
Near $(P,Q)$, we may express the $u=u_0$ contour 
\begin{equation}\label{u0contour}
2u_0 = h(X) - k(Y)
\end{equation}
as $Y(X)$ or as $X(Y)$.  
Differentiating (\ref{u0contour}) implicitly with respect to $X$ gives
$Y'(X) = h'(X)/k'(Y)$, and with respect to $Y$ gives $X'(Y) = k'(Y)/h'(X)$.
The tangent line to (\ref{u0contour})  at $(P,Q)$ is in a characteristic direction
if one of these derivatives is zero at $(P,Q)$.

\begin{ruul}\label{rule:tangentcontours}
Suppose $h$ and $k$ are differentiable, except possibly at points where the function has a vertical tangent line.
Then at all points $(P,Q)$ where $h'(P)/k'(Q)$ is zero or infinite,  the constant $u$ and $v$ contours of $\alpha$ given by $(U,V) = (h(X), k(Y))$ are tangent.
If $h'(P)/k'(Q)$ is real and nonzero, then the $u$ and $v$ contours through $(P,Q)$ intersect transversely.
\end{ruul}

In Fig.~\ref{c5c6}~(a), for example, where $h(t) = (t-1)^3 + 1$, the $u$ and $v$ contours are tangent to each other along the lines $X = \pm 1$ and $Y = \pm 1$.  In Fig.~\ref{c5c6}~(b), with $h(t) = t^2$, the contours are tangent along the characteristic axes.  
�
\bigskip

\centerline{\includegraphics[width=6in]{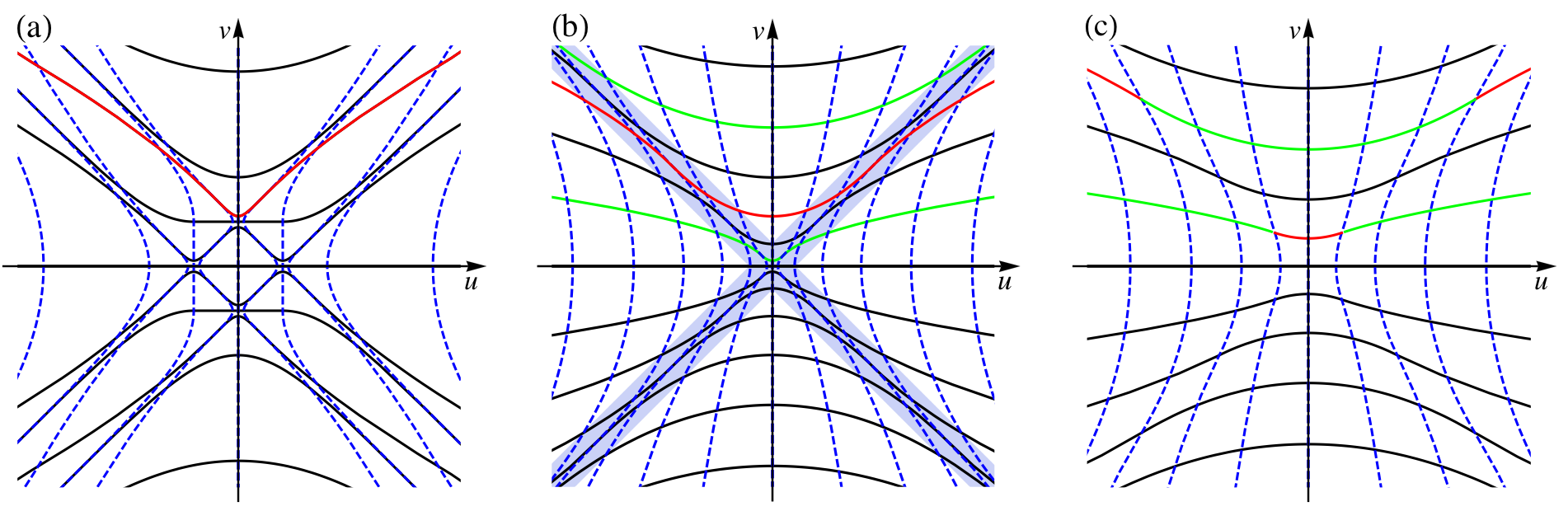}}

\begin{figure}[h]
\setlength{\abovecaptionskip}{2pt}
\caption{\small(a) Constant $x$ and $y$ contours of the Lorentz-conformal mapping  $(U,V)=(h(X), h(Y))$ with $h(t) = (t-1)^3 + 1$.  (b) Constant $x$ and $y$ contours of the unfolding from Construction~\ref{unfolding} with $h(t) = t^2$.  Contours are tangent along the characteristic axes.  (c) Constant $x$ and $y$ contours of $\alpha(c;X,Y)$ from Construction~\ref{cropping}, after cropping out the shaded region in (b), with $c=1$.  The coordinate lines now intersect transversely everywhere.  Red and green highlight how curves are re-joined after cropping.   
}
\label{c5c6}
\end{figure}

Given an increasing bijection $h$ on $[0, \infty)$ for which $h'(0)$ is zero or does not exist but $h'(t)$ is nonzero and finite for all $t\not=0$,
the mapping $\alpha$ from Construction~\ref{unfolding} is degenerate only along the characteristic axes, by Rule~\ref{rule:tangentcontours}.

\smallskip
{\bfseries Question:} 
How can this degeneracy be removed  to create a global Lorentz-conformal coordinate transformation that preserves the contours of $\alpha$ where they are transverse?
\smallskip

Construction~\ref{cropping} crops out strips around the characteristic axes and translates the remaining four wedges back to the origin so that the contours are re-joined along cropping lines.   
The width of the strip is determined by a parameter~$c$, which effects a homotopy of the original mapping by a family of Lorentz-conformal transformations for which the contours intersect transversely everywhere.  Fig.~\ref{c5c6}~(b,c) illustrates Construction~\ref{cropping} with $h(t)=t^2$ and $c=1$.  

\begin{const}[Cropping out degeneracies along characteristic axes]\label{cropping}
Let $h:\R\to\R$ be an odd bijection with $h'(t)>0$ for $t\not=0$.  
For each $c\geq0$, define $H_c: \mathbb{R} \rightarrow \mathbb{R}$~by
\[
H_c(t) = \left\{ \begin{array}{ll} h(t+c)-h(c), & t \geq 0  \\ 
h(t-c)-h(-c), & t \leq 0  . \end{array} \right.
\]
Then $\tilde\alpha : [0, \infty) \times \mathbb{R}^{1,1} \rightarrow \mathbb{R}^{1,1}$ defined by
\begin{equation}\label{alphahomotopy} 
\tilde\alpha(c;X,Y) = (H_c(X), H_c(Y))
\end{equation}
effects a homotopy of invertible Lorentz-conformal transformations $\alpha_c$ of the plane defined by $(U,V)=\tilde\alpha(c;X,Y)$, and for $c>0$, $\alpha_c$ is nondegenerate in that its constant $u$ and $v$ contours intersect transversely everywhere.
\end{const}

\subsection{$D_4$ symmetries and more unfoldings}\label{D4Symmetries}

The contour plots in Section~\ref{unfoldingsection} display evident symmetries.
For the 4:1 mapping in Fig.~\ref{unfoldingquadratic}(a), reflection about a characteristic axis or rotation by $\pi$ takes every contour onto itself, while reflection about a standard ($x$ or $y$) axis or rotation by $\pi/2$ or $3\pi/2$ preserves every contant-$v$ contour and interchanges the contours $u = u_0$ and $u = -u_0$.  
After unfolding, in Fig.~\ref{unfoldingquadratic}(b), the effects of these linear mappings on the contours are different.
The four reflections about the standard and characteristic axes, together with the powers of rotation by $\pi/2$, constitute the dihedral group $D_4$, of order 8. 
Motivated by these observations, we undertake a systematic study of the $D_4$ symmetries of Lorentz-conformal mappings. 

Consider the group $D_4$ of linear transformations of the $(X,Y)$-plane generated by reflection about the line $X=0$, which we denote by $\sigma$, and reflection about the line $x=0$, which we denote by $\tau$.  Then $\sigma \tau$ is rotation by $\pi/2$, and
$$D_4 = <\sigma, \tau \ | \ \sigma^2 = \tau^2 = (\sigma \tau)^4 = e >.$$ 
$D_4$ is the semidirect product of the (normal) rotation subgroup $T$, generated by $\sigma \tau$ and isomorphic to $\mathbb{Z}_4$, and the reflection subgroup $R_X$, generated by $\sigma$ and isomorphic to $\mathbb{Z}_2$.  We express this in tabular form as
\begin{equation}  \label{D4matrix}
D_4 = \mathbb{Z}_4 \rtimes \mathbb{Z}_2 = 
\left[ \begin{array}{rrrr}
e & \sigma \tau & (\sigma \tau)^2 & (\sigma \tau)^3 \\
\sigma &  \sigma (\sigma \tau) & \sigma (\sigma \tau)^2 & \sigma (\sigma \tau)^3 
\end{array} \right] = 
\left[ \begin{array}{llll}
e & \sigma \tau & \sigma \tau \sigma \tau & \tau \sigma \\
\sigma & \tau & \tau \sigma \tau & \sigma \tau \sigma 
\end{array} \right].
\end{equation}
Table~\ref{D4action} gives the action of $D_4$ on the plane in characteristic and standard coordinates.  

\bigskip 
\begin{table}[h] 
    \begin{center}
       \begin{tabular}{|l|l|l|l|} \hline 
 \multirow{2}{*}{$g \in D_4$} &  Angle of rotation  &  Image of  &  Image of  \\
  & or line of reflection & $(X,Y)$ & $(x,y)$ \\ \hline  \hline
$e$ & identity & $(X,Y)$ & $(x,y)$ \\ \hline
$\sigma \tau$ & 
rotation by ${\pi / 2}$ & $(-Y,X)$ & $(-y,x)$ \\ \hline
$\sigma \tau \sigma \tau = (\sigma \tau)^2$ &  
rotation by $\pi$ &  $(-X,-Y)$   & $(-x,-y)$  \\ \hline
$\tau \sigma = (\sigma \tau)^3$ &  
rotation by ${3\pi / 2}$ & $(Y,-X)$ & $(y,-x)$ \\ \hline
$\sigma$ & 
reflection across $X=0$ & $(-X,Y)$ & $(-y,-x)$  \\ \hline
$\tau = \sigma (\sigma \tau)$ & 
reflection across $x=0$ & $(Y,X)$ & $(-x,y)$  \\ \hline
$\tau \sigma \tau = \sigma (\sigma \tau)^2$  & 
reflection across $Y=0$ & $(X,-Y)$ & $(y,x)$  \\ \hline
$\sigma \tau \sigma = \sigma (\sigma \tau)^3$ & 
reflection across $y=0$ & $(-Y,-X)$ & $(x,-y)$ \\ \hline
\end{tabular}
    \end{center}
\caption{\small Action of $D_4$ on the $(x,y)$-plane in standard and characteristic coordinates.  For the action of $D_4$ on the $(u,v)$-plane, replace $X, Y, x, y$ by $U, V, u, v$.} \label{D4action}
\end{table}

Given a Lorentz-conformal mapping $\alpha$ and a transformation
$g \in D_4$ on the $(x,y)$-plane, we may ask what effect the action of $g$ has on the $u$ and $v$ contours of $\alpha$.  That is, 
 we seek all pairs $(g, g') \in D_4 \times D_4$ such that $\alpha \circ g = g' \circ \alpha $, where $g$ acts on the $(x,y)$-plane and $g'$ acts on the $(u,v)$-plane.
 We denote the set of all such pairs $(g, g')$ as $S_{\alpha}$:
$$
S_{\alpha} = \{(g, g') \in D_4 \times D_4 : \alpha \circ g = g' \circ \alpha \}.
$$
$S_{\alpha}$ is a subgroup of $D_4 \times D_4$; we refer to $S_{\alpha}$ as the {\it (full) symmetry group of $\alpha$}.  Any subgroup $S \subseteq S_{\alpha}$ is a {\it (partial) symmetry group of $\alpha$}.   We omit the words ``full" and ``partial" when the distinction is clear or not needed.  

For any full or partial symmetry group $S$ for some $\alpha$, the sets
$$S_1= \{g : (g,g') \in S \ \text{for some} \  g' \in D_4\}, \ \ \ 
S_{2} = \{g' : (g,g') \in S \ \text{for some} \  g \in D_4\}$$
are subgroups of $D_4$. 
Since the action of $D_4$ on the plane has trivial stablizer on every point that is not on a standard or characteristic axis, it follows that 
for every $g \in S_1$, there is a unique $g' \in S_2$ such that $(g, g') \in S$.  We may therefore consider the function
$\Phi : S_1  \rightarrow  S_2$, where $\Phi (g) = g'$, and verify that $\Phi $ is a surjective homomorphism.  It follows that
knowing the symmetry group $S$ is equivalent to knowing the homomorphism
$\Phi: S_1  \rightarrow  S_2$.  
For the full symmetry group $S_{\alpha}$, we refer to $\Phi_{\alpha} : S_{\alpha 1}  \rightarrow  S_{\alpha 2}$ as the {\it (full) symmetry homomorphism of $\alpha$}.  For a subgroup $S \subseteq S_{\alpha}$, we call $\Phi: S_1  \rightarrow  S_2$ a
{\it (partial) symmetry homomorphism of $\alpha$}.  
We may refer to a symmetry homomorphism $\Phi: S_1  \rightarrow  S_2$ as an ``$S_1$-symmetry."

Take, for example, a 4:1 mapping $\beta$: $(U,V) = (h(X), h(Y))$, where $h$ is even.  From Table~\ref{D4action} and the form of $\beta$, we see that $g = \sigma \tau$ takes $(h(X),h(Y))$ to $(h(-Y),h(X)) = (h(Y),h(X))$, so that $(U,V) \mapsto (V,U)$, or 
$(u,v) \mapsto (-u,v)$, which gives $g' = \tau$.  
Working this out for all $g \in D_4$, one
finds that the full symmetry group of $\beta$ is given by the homomorphism $D_4 \to \{ e, \tau\}$, where $\sigma \tau \mapsto \tau$ and $\sigma \mapsto e$.  One finds similarly that for the unfolding of $\beta$ via Construction~\ref{unfolding}, where $h=k$ and $h$ is odd, the symmetry group is the diagonal of $D_4 \times D_4$, given by the identity isomorphism $D_4 \to D_4$.  

\smallskip
{\bfseries Question:} Which homomorphisms $\Phi : S_1 \rightarrow  S_2$, of subgroups $S_1$ and $S_2$ of $D_4$, occur as symmetries of Lorentz-conformal mappings, and what are the conditions on the mapping imposed by each admissible symmetry?
\smallskip

For Lorentz-conformal mappings of the form $(U,V) = (h(X), k(Y))$, the answer is contained in Tables~\ref{SymTable1}, \ref{SymTable2}, and \ref{SymTable3}, which we now derive.  For Lorentz-conformal mappings of the form $(U,V) = (k(Y), h(X))$, the roles of $U$ and $V$ are interchanged. 
  
\begin{table}[h] 
    \begin{center}
       \begin{tabular}{|c|c|c|l|} \hline 
Subgroup & Generating set &  Isomorphic to  & Description \\ \hline  \hline
$D_4$  &   $\sigma, \, \sigma \tau$ 
& $\mathbb{Z}_4 \rtimes \mathbb{Z}_2$  &  Dihedral group of order 8 \\ \hline
$R_{XY}$ &    $\sigma, \, \tau \sigma \tau$
& $\mathbb{Z}_2 \times \mathbb{Z}_2$  & Reflections about characteristic axes\\ \hline
$R_{xy}$ &   $\tau, \, \sigma \tau\sigma$  
&   $\mathbb{Z}_2 \times \mathbb{Z}_2$ & Reflections about standard axes\\ \hline
 $T$ &   $\sigma \tau$ 
&  $\mathbb{Z}_4$  & Rotation subgroup \\ \hline
$R_X$  & $\sigma$  
&  $\mathbb{Z}_2$  & Reflections about $X = 0$ \\ \hline
$R_Y$  &  $\tau \sigma \tau$  
&  $\mathbb{Z}_2$  & Reflections about $Y = 0$ \\ \hline
$R_x$  &  $\tau$  
&  $\mathbb{Z}_2$  & Reflections about $x = 0$ \\ \hline
$R_y$  & $\sigma \tau \sigma$  
&  $\mathbb{Z}_2$  & Reflections about $y = 0$ \\ \hline
 $T_{\pi}$  &   $(\sigma \tau)^2$  
& $\mathbb{Z}_2$  & Rotations by $\pi$ \\ \hline
$e$ & $e$ & $I$ & Identity \\ \hline
\end{tabular}
    \end{center}
\caption{\small The subgroups of $D_4$, with generators and isomorphism type.  The notation corresponds to the action of the subgroup on the $(x,y)$-plane;  when the subgroup is regarded as acting on the $(u,v)$-plane, replace $X, Y, x, y$ by $U, V, u, v$.} \label{D4subgroups}
\end{table}

We ask which pairs $(g,g') \in D_4 \times D_4$ belong to symmetry groups of Lorentz-conformal mappings $\alpha$: $(U,V) = (h(X), k(Y))$.
With respect to the matrix of $D_4$ in (\ref{D4matrixagain}), the arrays (\ref{alphag}) and (\ref{galpha}) give the $(U,V)$ coordinates of $\alpha \circ g$ and $g' \circ \alpha$ at $(X,Y)$.

\begin{equation}  \label{D4matrixagain}
D_4 =  [g_{ij}] = 
\left[ \begin{array}{llll}
e & \sigma \tau & (\sigma \tau)^2 & (\sigma \tau)^3 \\
\sigma & \tau & \tau \sigma \tau & \sigma \tau \sigma 
\end{array} \right]
\end{equation}

\begin{equation}\label{alphag}
[\alpha \circ g_{ij}|_{(X,Y)}] = 
\left[\hspace{-4pt} \begin{array}{cccc}
\big(h(X),k(Y)\big) & \big(h(-Y),k(X)\big)  & \big(h(-X),k(-Y)\big) & \big(h(Y),k(-X)\big) \\
\vspace{-2ex}\\
\big(h(-X),k(Y)\big) & \big(h(Y),k(X)\big)  & \big(h(X),k(-Y)\big) & \big(h(-Y),k(-X)\big) 
\end{array} \hspace{-4pt}\right]
\end{equation}

\begin{equation}\label{galpha}
[g'_{ij} \circ \alpha|_{(X,Y)}] = 
\left[\hspace{-4pt} \begin{array}{cccc}
\big(h(X),k(Y)\big) & \big(\!\!-\!k(Y),h(X)\big)\!  & \big(\!\!-\!h(X),-k(Y)\big) & \big(k(Y),-h(X)\big) \\
\vspace{-2ex}\\
\big(\!\!-\!h(X),k(Y)\big)\! & \big(k(Y),h(X)\big)  & \big(h(X),-k(Y)\big) & \!\!\big(\!\!-\!k(Y),-h(X)\big) 
\end{array} \hspace{-4pt}\right]
\end{equation}

\bigskip 

Arrays (\ref{alphag}) and (\ref{galpha}) allow us to read off all classes of
Lorentz-conformal transformations $\alpha$: $(U,V) = (h(X), k(Y))$ whose symmetry groups are subgroups of $D_4 \times D_4$ and write down the symmetry homomorphisms for each class.  To determine whether $(g,g') \in D_4 \times D_4$ belongs to the symmetry group of some $\alpha$, we equate the entry $\alpha \circ g_{ij}$ in (\ref{alphag}) with the entry $g'_{ij} \circ \alpha$ in (\ref{galpha}).
For example, $(\sigma, \sigma)$ belongs to the symmetry group of $\alpha$ if and only if $h(-X) = -h(X)$ so that $h$ is odd.   

We rule out pairs $(g,g')$ that force $h$ or $k$ to be a constant function,
such as $(e, \tau)$, which gives $h(X) = k(Y)$ for all $X$ and $Y$.  Ruling out such pairs enforces that $g$ and $g'$ be both in columns 1 or 3, or both in columns 2 or 4.  Tables~\ref{SymTable1}, \ref{SymTable2},  and ~\ref{SymTable3} list all classes of Lorentz-conformal transformations $\alpha$ of the form $(U,V) = (h(X), k(Y))$ that have symmetry groups in $D_4 \times D_4$, together with the symmetry homomorphisms that are equivalent to the given condition(s) on $\alpha$.  The final symmetry in each row is the full symmetry homomorphism; it is implied by each partial symmetry that precedes it.  Partial symmetries that are not equivalent to the full symmetry are not listed; they may be read off by restricting the listed homomorphisms to proper subgroups of the domains.  

\bigskip 
\begin{table}[h] 
    \begin{center}
       \begin{tabular}{|c||c|c|c|} \hline 
Condition on & Equivalent partial & Equivalent partial & Full  \\
$(h,k)$ & $T$-symmetry & $R_{xy}$-symmetry & $D_4$-symmetry \\ \hline \hline
$h=k$, & $T \to R_u$ & $R_{xy} \to R_u$ & $D_4 \to R_u$  \\
$h$ even & $\sigma \tau \mapsto \tau$ & $(\tau, \sigma \tau \sigma) \mapsto (\tau, \tau)$ &
$(\sigma, \sigma \tau) \mapsto (e, \tau)$  \\ \hline 
$h=k$, & $T \to T$ & $R_{xy} \to R_{uv}$ & $D_4 \to D_4$ \\
$h$ odd &  $\sigma \tau \mapsto \sigma \tau$ &  $(\tau, \sigma \tau \sigma) \mapsto (\tau, \sigma \tau \sigma)$ &
$(\sigma, \sigma \tau) \mapsto (\sigma, \sigma \tau)$ \\ \hline \hline
$h=-k$, & $T \to R_v$ & $R_{xy} \to R_v$ & $D_4 \to R_v$ \\
$h$ even & $\sigma \tau \mapsto \sigma \tau \sigma$ & $(\tau, \sigma \tau \sigma) \mapsto (\sigma \tau \sigma, \sigma \tau \sigma) $ &
$(\sigma, \sigma \tau) \mapsto (e, \sigma \tau \sigma)$  \\ \hline 
$h=-k$, & $T \to T$ & $R_{xy} \to R_{uv}$ & $D_4 \to D_4$  \\
$h$ odd & $\sigma \tau \mapsto (\sigma \tau)^3$ &  $(\tau, \sigma \tau \sigma) \mapsto (\sigma \tau \sigma, \tau)$ &  
$(\sigma, \sigma \tau) \mapsto (\sigma, (\sigma \tau)^3)$ \\ \hline
\end{tabular}
    \end{center}
\caption{\small Classes of Lorentz-conformal mappings $(U,V) = (h(X), k(Y))$ that admit a full symmetry $\Phi_{\alpha} : D_4 \rightarrow  S_2$.  For each symmetry homomorphism $\Phi: S_1 \rightarrow  S_2$, the images of the generators of $S_1$ are given.  Each row lists every partial symmetry that is equivalent to the full symmetry and to the condition on $\alpha$.} \label{SymTable1}
\end{table}

\bigskip 

In Table~\ref{SymTable1}, the bijective transformations in the second and fourth rows ($h$ odd) are unfoldings of the 4:1 mappings in the first and third rows ($h$ even), respectively.
In each case, there are two options for the unfolding, depending on whether one restricts the even function to the positive or negative real line before extending it antisymmetrically.

\begin{figure}[h]\centerline{\includegraphics[width=6.1in]{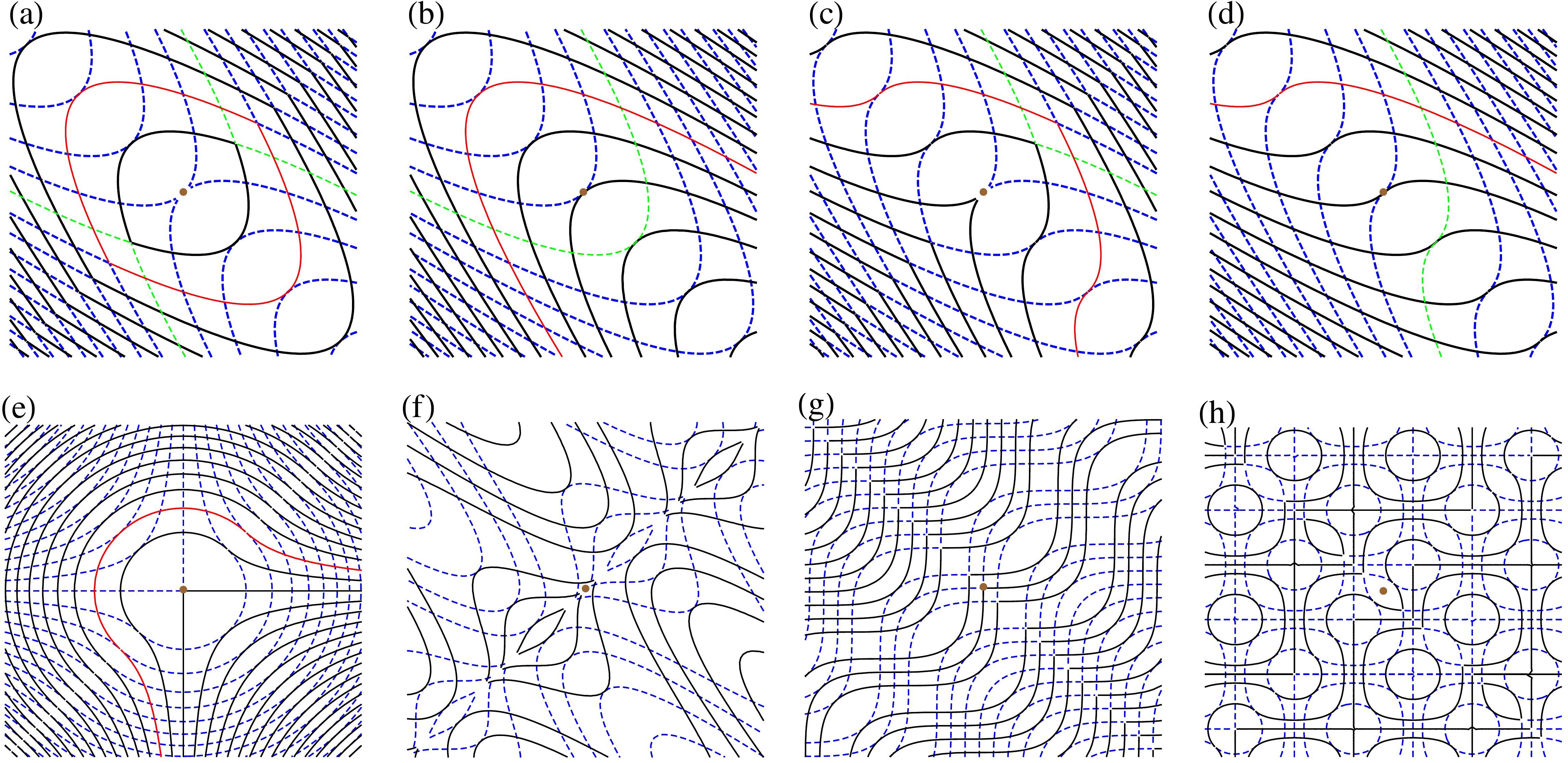}}
\caption{\small Contour plots of Lorentz-conformal transformations $(U,V) = (h(X), k(Y))$ illustrating the symmetry classes listed in Table~\ref{SymTable2}.  (a) $h(t)=t^2, k(t)=|t|$; (b) $h(t)=t^2, k(t)=t$; (c) $h(t)={\rm{sgn}}(t) t^2, k(t)=|t|$; (d) $h(t) = {\rm{sgn}}(t) t^2, k(t)=t$; 
(e) $h(t)=t^2, k(t) = {\rm{sgn}}(t) t^2$  (f) $h(t)= \sin(t),k(t)=\sqrt{|t|}$; (g) $h(t) = {\rm{sgn}}(t) \sin(t), k(t) = t$; (h) $h(t)={\rm{sgn}}(t)\sin(t), k(t) = {\rm{sgn}}(t) \cos(t)$.  In panels (a)-(d), one of the $u$-contours is marked in green, and in panels (a)-(e), one of the $v$-contours is marked in red.  The (brown) dot marks the origin.}
\label{table4figure}
\end{figure}

\begin{figure}[h]
\centerline{\includegraphics[width=6in]{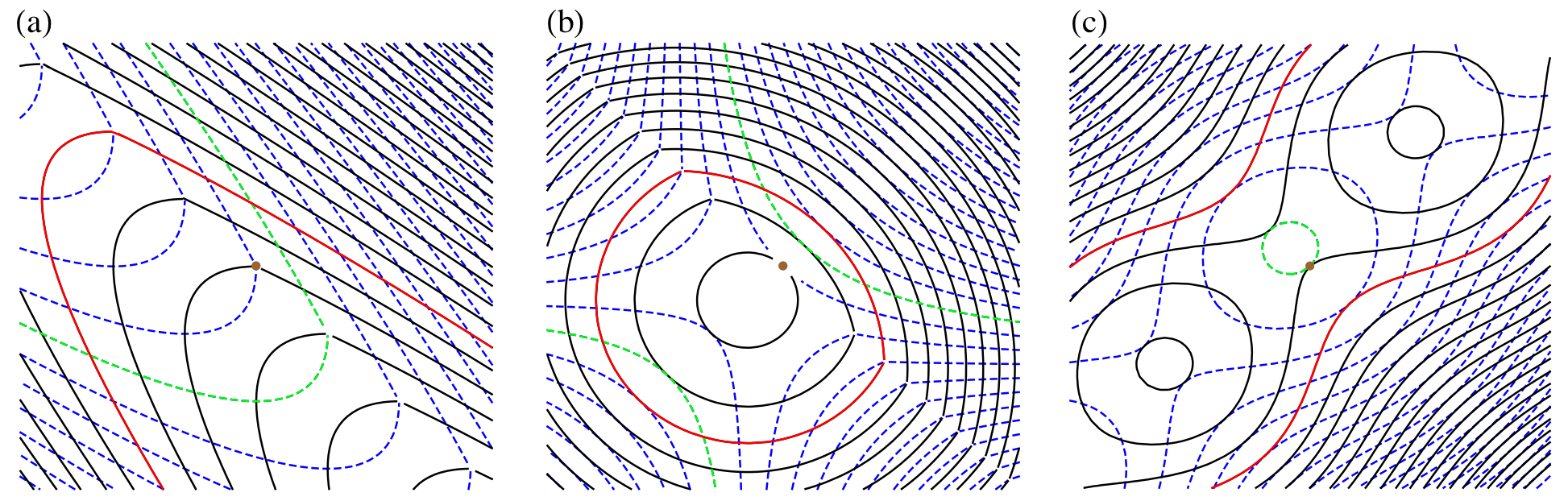}}
\caption{\small Contour plots of Lorentz-conformal transformations $\alpha$: $(U,V) = (h(X), k(Y))$ illustrating three of the symmetries in Table 5. (a) $h(t) = 2t+{\rm{sgn}}(t)t+t^2$, $k(t)=t$, (b)  $h(t) = 2t+{\rm{sgn}}(t)t+t^2$, $k(t)=t^2$, (c) $h(t) = \cos(t), k(t)=0.5(t^2-t)$.  In each panel, one of the $u$-contours is marked in green, and one of the $v$-contours is marked in red.  The (brown) dot marks the origin.}
\label{table5figure}
\end{figure}

Figs.~\ref{table4figure}~and~\ref{table5figure} respectively illustrate contour plots of transformations from the symmetry classes of Tables~\ref{SymTable2}~and~\ref{SymTable3}.  Compare Fig.~\ref{table4figure}~(e) with Fig.~\ref{unfoldingquadratic}(a,b).  Fig.~\ref{finalfigure} illustrates contour plots of Lorentz-conformal transformations that lack $D_4$-symmetries.

\begin{table}[h] 
    \begin{center}
	\begin{tabular}{|c||c|c|c|} \hline 
Condition on & Equivalent partial & Full  \\
$(h,k)$ & $T_{\pi}$--symmetry & $R_{XY}$--symmetry \\ \hline \hline
\multirow{2}{*}{$h$, $k$ even} & $T_{\pi} \to I$ & $R_{XY} \to I$  \\
  & $(\sigma \tau)^2 \mapsto e$ & $(\sigma, \tau \sigma \tau) \mapsto (e, e)$  \\ \hline 
\multirow{2}{*}{$h$ even, $k$ odd} & $T_{\pi} \to R_V$ & $R_{XY} \to R_V$  \\
& $(\sigma \tau)^2 \mapsto \tau \sigma \tau$ & $(\sigma, \tau \sigma \tau) \mapsto (e, \tau \sigma \tau)$  \\ \hline 
\multirow{2}{*}{$h$ odd, $k$ even} & $T_{\pi} \to R_U$ & $R_{XY} \to R_U$  \\
& $(\sigma \tau)^2 \mapsto \sigma$ & $(\sigma, \tau \sigma \tau) \mapsto (\sigma, e)$  \\ \hline 
\multirow{2}{*}{$h$, $k$ odd} & $T_{\pi} \to T_{\pi}$ & $R_{XY} \to R_{UV}$  \\
  & $(\sigma \tau)^2 \mapsto (\sigma \tau)^2$ 
  & $(\sigma, \tau \sigma \tau) \mapsto (\sigma, \tau \sigma \tau)$  \\ \hline 
\end{tabular}
 \end{center}
\caption{\small Classes of Lorentz-conformal transformations $(U,V) = (h(X), k(Y))$ with a full symmetry $\Phi : R_{XY} \rightarrow  I$.  The partial $T_{\pi}$-symmetry is equivalent to the full symmetry.}
\label{SymTable2}
\end{table}

The mappings in the second and third rows of Table~\ref{SymTable2} give 2:1 unfoldings of the 4:1 mappings ($h,k$ even) in the first row.  The invertible transformations in the fourth row ($h,k$ odd) are unfoldings of the non-invertible mappings in the first three rows.  


\bigskip 
\begin{table}[h] 
    \begin{center}
	\begin{tabular}{|c|c||c|c|} \hline 
Condition on & Full & Condition on & Full  \\
$(h,k)$ & symmetry & $(h,k)$ & symmetry\\ \hline \hline
\multirow{2}{*}{$h$ even} &  $R_X \to I$  & \multirow{2}{*}{$h(-t) = k(t)$}   &  $R_y \to R_u$  \\ 
 & $\sigma \mapsto e$ & & $\sigma \tau \sigma \mapsto \tau$ \\ \hline
\multirow{2}{*}{$h$ odd} &  $R_X \to R_U$  & \multirow{2}{*}{$h(-t) = -k(t)$}   &  $R_y \to R_v$  \\ 
 & $\sigma \mapsto \sigma$ & & $\sigma \tau \sigma \mapsto \sigma \tau \sigma$ \\ \hline
\hline
\multirow{2}{*}{$k$ even} &  $R_Y \to I$  & \multirow{2}{*}{$h=k$}   &  $R_x \to R_u$  \\ 
 & $\tau \sigma \tau \mapsto e$ & & $\tau \mapsto \tau$ \\ \hline
\multirow{2}{*}{$k$ odd} &  $R_Y \to R_V$  & \multirow{2}{*}{$h=-k$}   &  $R_x \to R_v$  \\ 
 & $\tau \sigma \tau \mapsto \tau \sigma \tau$ & & $\tau \mapsto \sigma \tau \sigma$ \\ \hline
\end{tabular}
 \end{center}
\caption{\small Lorentz-conformal transformations $(U,V) = (h(X), k(Y))$ for which $S_1$ has order 2 in the full symmetry homomorphism $\Phi: S_1 \rightarrow  S_2$.}
\label{SymTable3}
\end{table}

In Table~\ref{SymTable3}, the mappings with $h$ odd are unfoldings of those with $h$ even, and similarly for $k$.  If $h$ and $k$ are both increasing on $[0, \infty)$, then the mappings that satisfy
$h(-t) = k(t)$ are 4:1 and have unfoldings of the form $h(-t) = -k(t)$.  If one of $h$ and $k$ is increasing on $[0, \infty)$ and the other is decreasing, then the mappings that satisfy
$h(-t) = -k(t)$ are 4:1 and have unfoldings of the form $h(-t) = k(t)$.


\begin{center}
\includegraphics[width=6in]{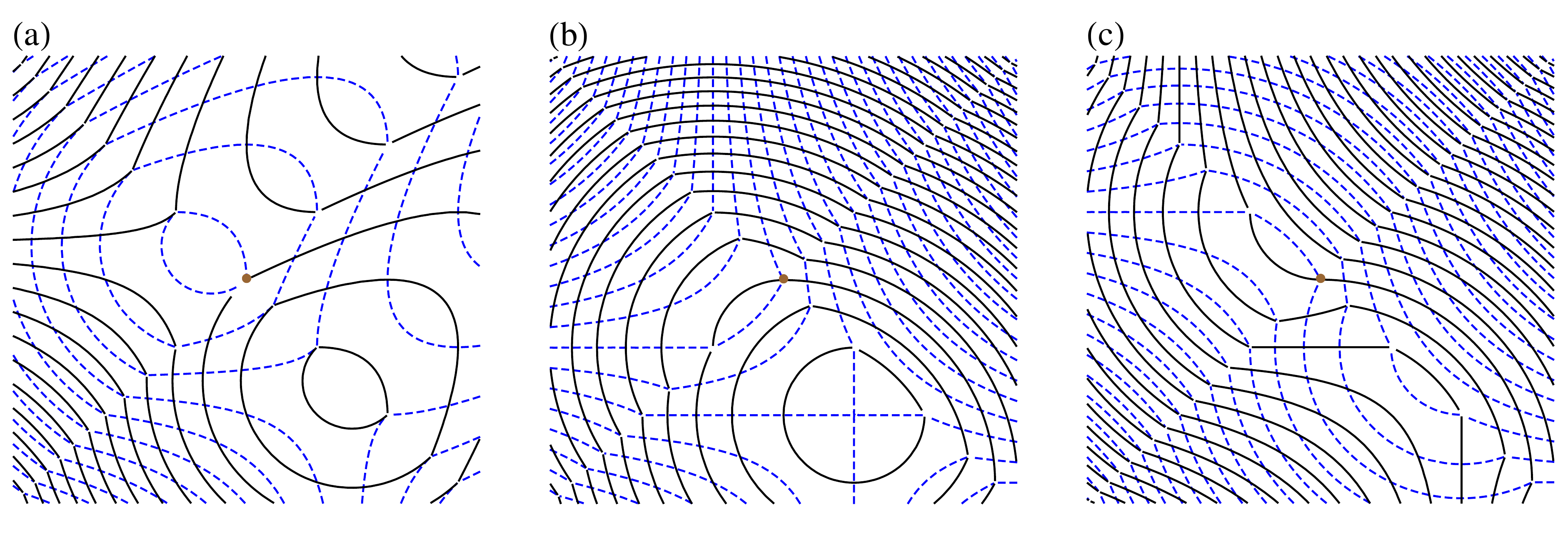}
\setlength{\abovecaptionskip}{2pt}
\captionof{figure}{\small Contour plots of Lorentz-conformal mappings $\alpha$: $(U,V) = (h(X), k(Y))$ that lack $D_4$ symmetries. (a) $h(t)=t$ for $t<0$ and  $t^2$ for $t>0$, $k(t) = 2t - {\rm{sgn}}(t)t-{\rm{sgn}}(t)t^2$. (b) $h(t)=2t+{\rm{sgn}}(t)t+t^2, k(t) = 2t - {\rm{sgn}}(t)t+t^2$  (compare with Fig.~\ref{Section3IntroductionFigure}(c)). (c)  $h(t)=2t+{\rm{sgn}}(t)t+ {\rm{sgn}}(t)t^2, k(t) = 2t - {\rm{sgn}}(t)t-{\rm{sgn}}(t)t^2$.  The (brown) dot marks the origin. }
\label{finalfigure}
\end{center}

\end{document}